\documentclass[11pt]{article}
\usepackage{amsmath,amsthm,amssymb}
\usepackage{caption,graphicx,subfig}
\usepackage[margin=1in]{geometry}
\usepackage{paralist,array}
\usepackage[sort&compress,numbers]{natbib}
\usepackage{setspace}

\setcaptionmargin{0.25in}

\makeatletter\@addtoreset{equation}{section}\makeatother

\theoremstyle{definition}

\newcolumntype{M}[1]{>{\centering\arraybackslash}m{#1}}

\newcommand{\drm}{\mathrm{d}}
\newcommand{\R}{\mathbb{R}}
\newcommand{\xv}{{\bf x}}
\newcommand{\Xv}{{\bf X}}
\newcommand{\Av}{{\bf A}}
\newcommand{\Bv}{{\bf B}}
\newcommand{\Cv}{{\bf C}}
\newcommand{\Iv}{{\bf I}}
\newcommand{\Jv}{{\bf J}}
\newcommand{\Kv}{{\bf K}}
\newcommand{\Mv}{{\bf M}}
\newcommand{\Pv}{{\bf P}}
\newcommand{\Qv}{{\bf Q}}
\newcommand{\Rv}{{\bf R}}
\newcommand{\Yv}{{\bf Y}}
\newcommand{\uv}{{\bf u}}
\newcommand{\Fv}{{\bf F}}

\newcommand{\Xiv}{{\bf \Xi}}
\newcommand{\xiv}{{\bf \xi}}
\newcommand{\Thetav}{{\bf \Theta}}
\newcommand{\thetav}{{\bf \theta}}

\newcommand*\samethanks[1][\value{footnote}]{\footnotemark[#1]}

\begin{document}

\title{\LARGE{\textbf{Data-Driven Stabilization of Periodic Orbits}}} 

\author{
Jason J. Bramburger\thanks{Department of Applied Mathematics, University of Washington, Seattle, WA, 98105}
\and
J. Nathan Kutz\samethanks  
\and 
Steven L. Brunton\thanks{Department of Mechanical Engineering, University of Washington, Seattle, WA, 98105}
}

\date{}
\maketitle

\begin{abstract}
Periodic orbits are among the simplest non-equilibrium solutions to dynamical systems, and they play a significant role in our modern understanding of the rich structures observed in many systems. For example, it is known that embedded within any chaotic attractor are infinitely many unstable periodic orbits (UPOs) and so a chaotic trajectory can be thought of as `jumping' from one UPO to another in a seemingly unpredictable manner. A number of studies have sought to exploit the existence of these UPOs to control a chaotic system. These methods rely on introducing small, precise parameter manipulations each time the trajectory crosses a transverse section to the flow. Typically these methods suffer from the fact that they require a precise description of the Poincar\'e mapping for the flow, which is a difficult task since there is no systematic way of producing such a mapping associated to a given system. Here we employ recent model discovery methods for producing accurate and parsimonious parameter-dependent Poincar\'e mappings to stabilize UPOs in nonlinear dynamical systems. Specifically, we use the sparse identification of nonlinear dynamics (SINDy) method to frame model discovery as a sparse regression problem, which can be implemented in a computationally efficient manner. This approach provides an explicit Poincar\'e mapping that faithfully describes the dynamics of the flow in the Poincar\'e section and can be used to identify UPOs. For each UPO, we then determine the parameter manipulations that stabilize this orbit. The utility of these methods are demonstrated on a variety of differential equations, including the R\"ossler system in a chaotic parameter regime.                
\end{abstract}


\section{Introduction}

Since their inception by Henri Poincar\'e at the turn of the twentieth century, return maps, or Poincar\'e maps as they are now commonly referred, have significantly influenced our understanding of recurrent and chaotic dynamical systems. Poincar\'e showed that one can understand the dynamics of a system in phase space not by looking at the full trajectory, but by tracking where trajectories intersect a lower-dimensional subspace, transverse to the flow. Such a subspace is now referred to as a Poincar\'e section, with the mapping that iterates between successive intersections of this section called the Poincar\'e mapping. Using a Poincar\'e mapping has the effect that continuous time-dependent behaviour of a dynamical system is translated into a discrete iterative process with a lower-dimensional phase-space. Importantly, fixed points and periodic/cyclic orbits of this mapping manifest themselves as periodic orbits of the continuous-time dynamical system, and therefore Poincar\'e mappings provide an accessible method of understanding the flow on and near periodic orbits in phase space, many of which are unstable and ultimately responsible for the rich, dynamical structures manifest in many systems.
Indeed, as we demonstrate, by exploiting data-driven methods for identifying such unstable periodic orbits, we can then stabilize a system by prescribing appropriate perturbations to its parameters to achieve a desired dynamic outcome. 

It is through discrete dynamical systems, such as Poincar\'e mappings, that much of our modern understanding of chaos comes from as well~\cite{Devaney,OttBook}. Concepts such as horseshoe mappings and symbolic dynamics now inform a significant portion of our comprehension of the geometry of a chaotic attractor. Particularly, it is known that any chaotic attractor contains an infinite number of unstable periodic orbits (UPOs) and that many dynamical averages such as Lyapunov exponents, entropy, and fractal dimensions can be expressed in terms of a weighted sum over the embedded UPOs~\cite{Artuso,Auerbach,Cvitanovic,OttBook}. 
UPOs play an important role in our modern understanding of several complex dynamical systems, such as turbulent fluid flow and transport in the solar system.  
In turbulence, UPOs guide the spatio-temporal chaos observed in spatially-extended fluid models such as the Navier--Stokes equations~\cite{Budanur,Cvitanovic2,Fazendeiro,Franceschini,Lucas,petrov1996nonlinear,eckhardt2007turbulence,Yalniz,suri2020capturing,graham2020exact}. 
There have also been considerable advances in energy efficient space mission design that take advantage of UPOs in the solar system~\cite{koon2000heteroclinic,gomez2004connecting,Dellnitz:2005}. 

In their seminal work, Ott, Grebogi, and Yorke demonstrated that these UPOs could be leveraged to control chaotic trajectories through small, precise parameter manipulations at each iterate of the mapping~\cite{Ott}. They proposed that when a chaotic trajectory comes sufficiently close to an element of a desired UPO, well-chosen parameter perturbations could be applied to keep the trajectory close to the chosen UPO, thus stabilizing the unstable orbit. They further commented that the ergodicity of orbits on the chaotic attractor would guarantee that eventually any chaotic trajectory comes sufficiently close to any element of the desired UPO, upon which the control algorithm can then be implemented. Hence, their method constitutes {\em closed-loop} control since the action from the controller is dependent on measurements of the system at each iteration. 

Since the original work of~\cite{Ott}, a number of advances have been made to simplify and extend algorithms for controlling chaos. From these extensions, the control of chaos has been successfully applied to cardiac rhythms~\cite{Weiss}, mechanical systems~\cite{Fradkov}, synchronizing electrical circuits~\cite{Liao}, satellite systems~\cite{Farid}, and more. The reviews~\cite{Andrievskii,Andrievskii2,Boccaletti,Fradkov2} and the references therein provide an overview of the many advances  in this area, as well as a more complete discussion of the applications. 

Despite the success of controlling chaos through these methods, there are still significant limitations that must be addressed before they can be applied to a broader range of chaotic systems. Primary among these issues is that of identifying an explicit Poincar\'e mapping. That is, aside from a few simple examples, there remains the challenge of obtaining an explicit parameter-dependent mapping that governs the iterates of a system inside a Poincar\'e section -- even for many simple dynamical systems. As the numerous studies on chaos control have shown, a closed-form Poincar\'e mapping for a given system can be used to determine the location of UPOs and the exact parameter perturbations necessary for control. Therefore, the most crucial piece to controlling chaos in the vein of~\cite{Ott} is identifying the Poincar\'e mapping associated to the given dynamical system.    

In this paper, we demonstrate that a recent method for data-driven discovery of Poincar\'e maps~\cite{Bramburger} can be used to overcome this barrier to stabilizing UPOs of nonlinear dynamical systems. We employ the {\em sparse identification of nonlinear dynamics} (SINDy) method~\cite{SINDy} for the discovery of the mappings, which has been proven to be a robust model-discovery algorithm and has been applied to discovering not only Poincar\'e mappings, but multiscale dynamics~\cite{Bramburger2,Champion}, partial differential equations~\cite{Rudy}, boundary value problems~\cite{BVP}, equations of rational expressions~\cite{Kaheman}, and conservation laws~\cite{Kaiser}. In the present application, we demonstrate that the SINDy method provides an explicit mapping that governs the dynamics of iterates within the Poincar\'e section, for which UPOs manifest themselves as unstable fixed points or cyclic orbits of this mapping. With a closed-form Poincar\'e mapping, UPOs can then be identified using simple root-finding techniques. Furthermore, the discovered mapping not only provides the location of UPOs, but can also be used to determine approximate stable and unstable directions associated to these periodic orbits. We will see in the many examples provided herein that applying SINDy to training data gathered in a neighbourhood of a focal parameter value leads to a faithful representation of the true Poincar\'e section dynamics that can be employed to both find and stabilize UPOs in a computationally efficient manner.           

\begin{figure} 
\center
\includegraphics[width = \textwidth]{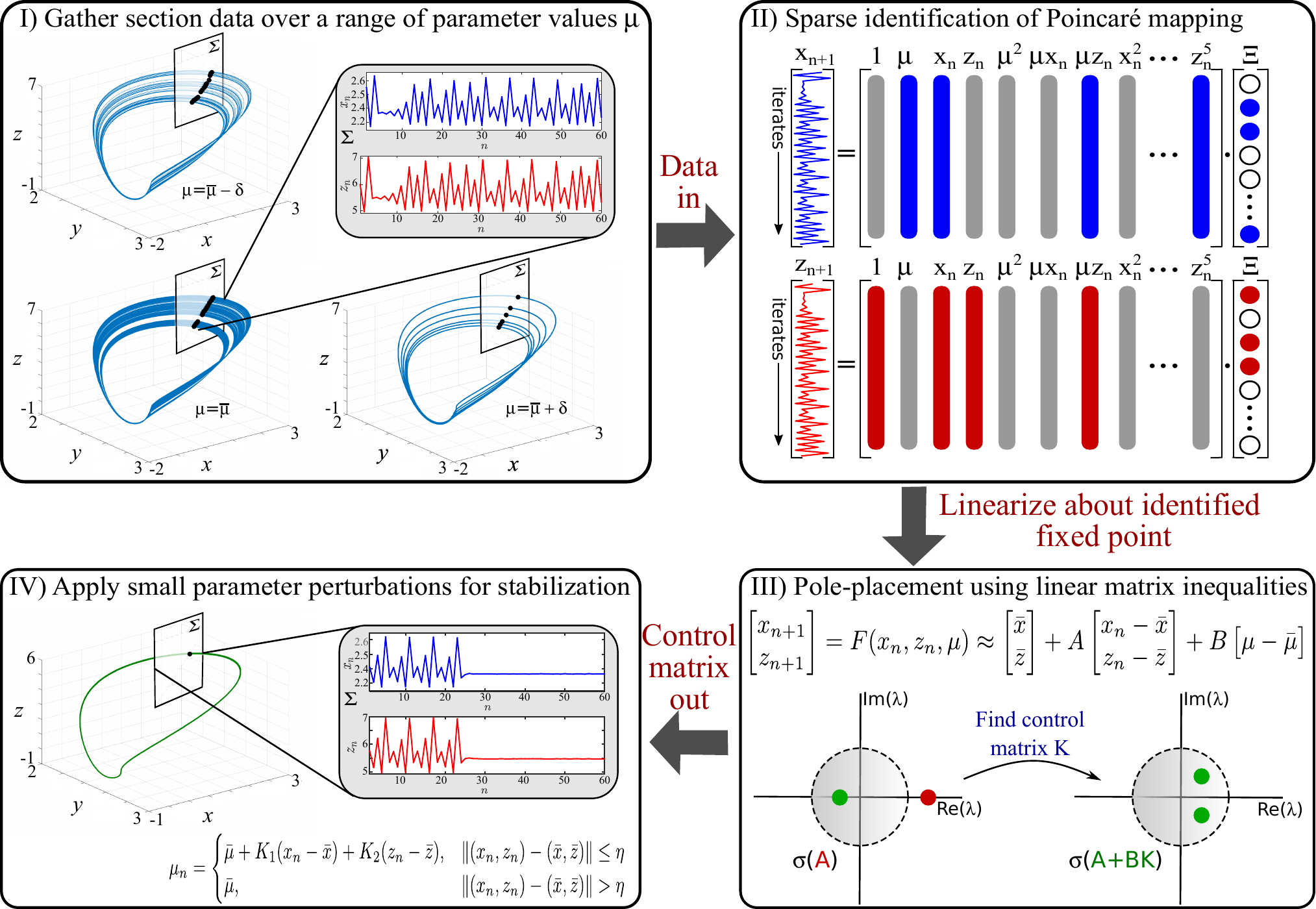} 
\caption{A visual overview of the stabilization procedure presented herein applied to a period 1 orbit. Given a parameter-dependent ODE, we first collect iterates in the Poincar\'e section, $\Sigma$, here a subset of the plane $y = 0$, over a range of parameter values near the focal value $\mu = \bar{\mu}$. {This provides the section data} $\xv_n = [x_n,z_n]$ {in the two remaining scalar variables $x$ and $z$.} Using the section data we then apply the SINDy method to discover a parsimonious parameter-dependent Poincar\'e mapping $[x_{n+1},z_{n+1}]^T = F(x_n,z_n,\mu)$, which can be analyzed explicitly. Upon identifying an unstable fixed point, $(\bar x, \bar z)$, of the discovered Poincar\'e mapping, we use linear matrix inequalities to find a control matrix $\Kv$ that stabilizes this unstable orbit through small perturbations of the parameter at each map iteration. These small perturbations can then be applied to stabilize UPOs in the original ODE by adjusting the parameter each time the trajectory intersects the Poincar\'e section near the UPO.}
\label{fig:Algorithm}
\end{figure} 

In Figure~\ref{fig:Algorithm}, we provide a visual illustration of the stabilization procedure, which is the focus of this article, applied to a fixed point of the Poincar\'e mapping. The method requires one to define a Poincar\'e section, $\Sigma$, that is transverse to the flow of a given nonlinear dynamical system and then obtain Poincar\'e section data, {written compactly as the vector iterates} ${\xv}_n$, over a range of parameter values about a focal value, say $\bar\mu$. This gathered data represents the training data required to seed the SINDy method for the discovery of a parameter-dependent discrete dynamical system, $\xv_{n+1} = \Fv(\xv_n,\mu)$, that maps iterates from one point in the Poincar\'e section to the next. 

Supposing that we have identified a fixed point $\bar\xv$ at $\mu = \bar\mu$, which corresponds to a UPO of the continuous-time dynamical system that we wish to stabilize\footnote{The extension to stabilizing cyclic orbits of the Poincar\'e map is handled similarly later in the manuscript.}, we linearize the discovered mapping $\Fv$ about $(\xv,\mu) = (\bar\xv,\bar\mu)$, {to obtain the Jacobian matrices evaluated at $(\bar\xv,\bar\mu)$, $\Av$ and $\Bv$, resulting from differentiation of $\Fv$ with respect to $\xv$ and $\mu$, respectively}. Since $\bar\xv$ is taken to correspond to a UPO, it follows that at least one of the eigenvalues of $\Av$ lies outside the unit circle in the complex plane, thus implying that the fixed point $\bar\xv$ is linearly unstable. To stabilize this orbit, we apply the method of~\cite{Romeiras} by introducing small parameter manipulations at each intersection of the continuous-time trajectory with the Poincar\'e section. {That is, we wish to obtain a control matrix $\Kv$ so that the eigenvalues of $\Av + \Bv\Kv$ are entirely contained in the unit circle of the complex plane. Having found the matrix $\Kv$ we may then return to the original continuous-time dynamical system and apply the appropriate parameter manipulation each time the trajectory transversely intersects the Poincar\'e section.} {These parameter are described in detail in} Section~\ref{sec:LMI} {as well as} Figure~\ref{fig:Algorithm}. 

From the above description of the stabilization algorithm, we see that obtaining an appropriate control matrix $\Kv$ to stabilize the trajectory is an important aspect of the procedure that must be informed by the linearization of the discovered mapping about the iterates of the UPO in the section. We note that there are a number of efficient methods which could be employed to finding the control matrix $\Kv$~\cite{Dullerud,Ogata,Sontag,DDBook,sp:book,Astrom2010book}. Here we opt to follow the methods of Parrilo~\cite{Parrilo} to see how the control matrix $\Kv$ can be obtained in a principled manner using linear matrix inequalities (LMIs). LMIs have long been applied to control systems~\cite{Boyd} and have the advantage that they can be implemented numerically as semidefinite programs. This eases the analysis required by the user to stabilize UPOs since it does not require one to solve complicated nonlinear equations that may arise from attempting to find $\Kv$ by hand. {Recent similar investigations have used Taylor expansions to obtain approximations of Poincar\'e mappings, while also employing LMIs to stabilize fixed points in the section} \cite{Rev1,Rev2,Rev3,Rev4}. In this work we will extend the work of~\cite{Parrilo} to show that LMIs can be used to determine control matrices that stabilize cyclic orbits of the Poincar\'e mapping, going beyond just the original application to fixed points.  

The rest of this paper is organized as follows. In Section~\ref{sec:SINDy} we present the SINDy method for discovering parameter-dependent Poincar\'e mappings. The presentation is supplemented with a brief historical overview of the method and a discussion of why the SINDy method is an appropriate model-discovery algorithm for the task at hand. Then in Section~\ref{sec:LMI} we go over the LMI relaxation procedure of~\cite{Parrilo} for obtaining a control matrix $\Kv$. The original method for stabilizing fixed points of a mapping is presented in \S~\ref{subsec:FixedPt}, while our simple extension to cyclic orbits of a mapping is presented in \S~\ref{subsec:Cycle}. Following these theoretical presentations of the two major components of the method, in Section~\ref{sec:Examples} we present a number of applications. This includes demonstrating the utility of the LMI framework for stabilizing orbits of the chaotic H\'enon map and then moving to the full algorithm which includes discovery of a Poincar\'e mapping. Our applications of the full method features a toy model which exhibits an isolated UPO and two chaotic systems, including the R\"ossler system. Section~\ref{sec:ParamInd} briefly discusses how the method is adapted for parameter-independent models. The parameter-independent method is applied to a classical problem of having a satellite rest at an unstable Lagrange point of the Earth and its moon by applying discrete thrusts at regular temporal intervals. The paper concludes in Section~\ref{sec:Discussion} with a review of our findings and a discussion of some avenues for future work.


\section{Discovery of Parameter-Dependent Poincar\'e Maps}\label{sec:SINDy}

Given the flow $\xv(t)$ of a continuous-time dynamical system, Poincar\'e proposed defining a lower-dimensional hyperplane $S$ that is transverse to this flow to simply track the intersection of the flow with $S$. The result of this process is a discrete sequence
\begin{equation}
	\{\xv_n = \xv(t_n):\ \xv(t_n)\in S,\ t_n > t_{n-1}\}
\end{equation}
which comprise the successive iterates in the section. As discussed in the introduction, we are interested in parameter-dependent dynamics and therefore we will assume that such a sequence of iterates in a Poincar\'e section can be obtained over a range of parameter values, leading to the discrete sequences of data snapshots $\{\xv_1^{(j)},\dots,\xv_{N_j}^{(j)}\}\subset\R^d$ for the sampled parameter values $\mu_j \in\R^p$ with $j = 1,\dots,M$. Our goal is to determine a parameter-dependent mapping so that when $\mu = \mu_j$ we map vectors $\xv_n^{(j)}$ to $\xv_{n+1}^{(j)}$. Such a mapping is referred to as a Poincar\'e mapping since it maps successive iterates of the Poincar\'e section into each other. As noted in the introduction, explicit representations of Poincar\'e mappings are rarely obtained for systems of interest. Therefore, we will leverage the model discovery algorithm SINDy~\cite{Bramburger,SINDy} to discover these parameter-dependent mappings.

We may begin by augmenting the data to incorporate the parameter as an extra $p$-variables in the training data. For each $j = 1,\dots,M$ we consider the sequences $\{(\xv_1^{(j)},\mu_j),\dots,(\xv_{N_j}^{(j)},\mu_j)\}\subset\R^d\times\R^p$ with the assumption that $(\xv_n^{(j)},\mu_j) \mapsto (\xv_{n+1}^{(j)},\mu_j)$, i.e. the parameter dynamics are trivial. We refer to this mapping as $\Fv:\R^d \times \R^p \to\R^d \times \R^p$ throughout this section. To discover $\Fv$ we apply SINDy, which frames model discovery as a sparse regression problem. We first stack the data into matrices 
\begin{equation}
	\Xv_1 := \begin{bmatrix} \xv_1^{(1)} & \mu_1 \\ \vdots & \vdots \\ \xv_{N_1-1}^{(1)} & \mu_1 \\ \xv_1^{(2)} & \mu_2 \\ \vdots & \vdots \\ \xv_{N_2-1}^{(2)} & \mu_2 \\ \vdots & \vdots \\ \xv_1^{(M)} & \mu_M \\ \vdots & \vdots \\ \xv_{N_M-1}^{(M)} & \mu_M \end{bmatrix}, \quad \Xv_2 := \begin{bmatrix} \xv_2^{(1)} & \mu_1 \\ \vdots & \vdots \\ \xv_{N_1}^{(1)} & \mu_1 \\ \xv_2^{(2)} & \mu_2 \\ \vdots & \vdots \\ \xv_{N_2}^{(2)} & \mu_2 \\ \vdots & \vdots \\ \xv_2^{(M)} & \mu_M \\ \vdots & \vdots \\ \xv_{N_M}^{(M)} & \mu_M \end{bmatrix}
\end{equation} 
belonging to $\R^{(N-M) \times (d+p)}$, where $N = \sum_{j = 1}^M N_j$. From the assumptions on the data $(\xv_n^{(j)},\mu_j)$, by construction $\Fv$ maps each row of $\Xv_1$ into the same row of $\Xv_2$. Although $\Fv$ is unknown, we may construct a library of $q \geq 1$ candidate functions for which $\Fv$ belongs to the linear span of these functions. Let us denote these candidate model functions $\thetav_1,\dots,\thetav_q:\R^{d+p}\to\R^{d+p}$ so that the assumption $\Fv \in \mathrm{span}\{\thetav_1,\dots,\thetav_q\}$ {becomes}
\begin{equation}\label{Fspan}
	\Fv(\xv,\mu) = \xiv_1\thetav_1(\xv,\mu) + \dots + \xiv_q\thetav_q(\xv,\mu)
\end{equation}
{where the $\xiv_1,\dots,\xiv_q\in\R$ make up the coefficients in the linear span. For example, with $d = p = 1$, one could consider polynomial basis functions up to quadratic order in a single variable $x$ and a single parameter $\mu$, thus giving a library of the form}
\begin{equation}
\begin{bmatrix}
	1 \\ 0		
\end{bmatrix}, \begin{bmatrix}
	x \\ 0		
\end{bmatrix}, \begin{bmatrix}
	\mu \\ 0		
\end{bmatrix}, \begin{bmatrix}
	0 \\ \mu		
\end{bmatrix}, \begin{bmatrix}
	x^2 \\ 0		
\end{bmatrix}, \begin{bmatrix}
	\mu x \\ 0		
\end{bmatrix}, \begin{bmatrix}
	\mu^2 \\ 0		
\end{bmatrix},
\end{equation} 
{where we remind the reader that the parameter dynamics, i.e. the second component, are assumed to be trivial. Due to these trivial parameter dynamics, we do not include them in our description of the discovered mappings in the applications that follow.}

{Now, we can couple the assumption} \eqref{Fspan} {and the fact that rows of $\Xv_1$ are mapped by $\Fv$ into rows of $\Xv_2$. Take the first row as a specific example. Then, by assumption we have}
\begin{equation}
\begin{bmatrix}
	\xv_2 \\ \mu_1
\end{bmatrix} = \xiv_1\thetav_1(\xv_1,\mu_1) + \dots + \xiv_q\thetav_q(\xv_1,\mu_1),
\end{equation} 
{giving a linear equation for the coefficients $\Xiv = [\xiv_1\ \dots\ \xiv_q]^T \in \R^q$. We can repeat this process for each row of $\Xv_1$ and $\Xv_2$, resulting in the linear system}
\begin{equation}\label{SINDyMatrix}
	\Xv_2 = \Thetav(\Xv_1) \Xiv
\end{equation}
{where}
\begin{equation}
\begin{split}
		&\Thetav(\Xv_1) \\ &= \begin{bmatrix}
			\thetav_1(\xv_1^{(1)},\mu_1) & \thetav_2(\xv_1^{(1)},\mu_1) & \cdots & \thetav_q(\xv_1^{(1)},\mu_1) \\
			\vdots & \vdots & \ddots & \vdots \\
			\thetav_1(\xv_{N_1-1}^{(1)},\mu_1) & \thetav_2(\xv_{N_1-1}^{(1)},\mu_1) & \cdots & \thetav_q(\xv_{N_1-1}^{(1)},\mu_1) \\
			\thetav_1(\xv_1^{(2)},\mu_2) & \thetav_2(\xv_1^{(2)},\mu_2) & \cdots & \thetav_q(\xv_1^{(2)},\mu_2) \\
			\vdots & \vdots & \ddots & \vdots \\
			\thetav_1(\xv_{N_2-1}^{(2)},\mu_2) & \thetav_2(\xv_{N_2-1}^{(2)},\mu_2) & \cdots & \thetav_q(\xv_{N_2-1}^{(2)},\mu_2) \\
			\vdots & \vdots & \ddots & \vdots \\
			\thetav_1(\xv_{1}^{(M)},\mu_M) & \thetav_2(\xv_{1}^{(M)},\mu_M) & \cdots & \thetav_q(\xv_{1}^{(M)},\mu_M) \\
			\vdots & \vdots & \ddots & \vdots \\
			\thetav_1(\xv_{N_M-1}^{(M)},\mu_M) & \thetav_2(\xv_{N_M-1}^{(M)},\mu_M) & \cdots & \thetav_q(\xv_{N_M-1}^{(M)},\mu_M) \\
		\end{bmatrix}
\end{split}
\end{equation}   
{is a matrix of the library functions evaluated at the entries of $\Xv_1$}.

In practice we require a large number of measurements for each parameter value, thus giving that $N - M \gg 1$ and making \eqref{SINDyMatrix} an overdetermined linear system. The early works of Small et al.~\cite{Small} and Yao and Bollt~\cite{Yao} formulated model discovery in a similar manner resulting in models that potentially include all candidate functions after solving for $\Xiv$. The over-determined linear system \eqref{SINDyMatrix} has a number of regularized solutions, but in practice our goal is to have the mapping $\Fv$ written as the linear span of as few candidate functions as possible, so a sparseness requirement can be imposed on the coefficients of $\Xiv$. In the present context, this sparsity requirement is motivated by the following three points: 1) many physical systems can be written in a sparse basis of candidate functions, 2) using fewer candidate functions to describe $\Fv$ will make the discovered mapping easier to analyze in what follows, and 3) promoting sparsity will help to eradicate small terms in the solution of \eqref{SINDyMatrix} which could be attributed to numerical error in the training data.  

There are many variants of sparse regression that could be employed to solve \eqref{SINDyMatrix}, all of which attempt to approximate the solution of an NP-hard $\ell_0$-optimization process~\cite{Natarajan}. In this manuscript we employ sequential threshold least squares~\cite{SINDy} which is a proxy for $\ell_0$-optimization~\cite{Zheng}, has convergence guarantees~\cite{Zhang}, and out performs LASSO~\cite{Su,Tibshirani} in most cases at significantly reduced computational expense. To apply this method we begin with the least-squares solution of \eqref{SINDyMatrix}, i.e. minimizing $\|\Xv_2 - \Thetav(\Xv_1)\Xiv\|_2$ in $\Xiv = [\xiv_1\ \dots\ \xiv_q] \in \R^q$. {We then introduce a sparsity threshold, $\lambda > 0$ and define a new unknown matrix $\tilde{\Xiv} =  [\tilde\xiv_1\ \dots\ \tilde\xiv_q] \in \R^q$ with $\tilde{\xiv}_j = 0$ if the least-squares solution $\Xiv$ is such that $|\xiv_j|\leq \lambda$. This then moves one to another overdetermined linear system of equations}
\begin{equation}\label{LinSyst}
	\Xv_2 = \Thetav(\Xv_1) \tilde\Xiv
\end{equation} 
{which potentially has fewer degrees of freedom due to fixing some coefficients of $\tilde\Xiv$ to zero. We then repeat the procedure by obtain the least-squares solution of} \eqref{LinSyst} {with the dimensionally-reduced vector $\tilde{\Xiv}$.} This process is iterated with $\tilde{\Xiv}$ replacing $\Xiv$ until a sufficiently sparse model is obtained.  

When implementing the above method there are two major components that must be chosen prior to implementing the method: the candidate functions which make up the library and the sparsity parameter. The choice of these two pieces of the method are intimately related. For example, one can image a scenario where the true mapping contains a trigonometric term, while the library consists solely of monomials. It is not unreasonable to suppose that the discovered mapping would attempt to build a series representation of the trigonometric term and the choice of the sparsity parameter effectively tells the discovery process where to truncate this series. Of course, taking $\lambda$ too large will provide too few terms in the series, while taking $\lambda$ too small has the potential to allow numerical error in cultivating the training data to manifest itself as small coefficients on library functions that should not be present in the discovered mapping. Conversely, including the appropriate trigonometric term in the library of functions could potentially provide greater forgiveness when choosing a value for $\lambda$. Such a phenomenon was exemplified in~\cite[Section~IV C]{Bramburger2} where polynomial terms were used to approximate dynamics governed by a sine function. It can be a difficult task choosing the appropriate functions to include in the library of candidate functions and at the time of writing this manuscript there does not appear to be any prescriptions that are broadly applicable. Although this is a potential limitation of the method, we will see in Section~\ref{sec:Examples} that for the examples considered herein, using a simple library of monomials can reliably capture the qualitative dynamics of a desired Poincar\'e mapping.   

We conclude this section with the following remark. If there is a particular parameter value of interest, say $\bar{\mu}$, that one wishes to focus their analysis in a neighbourhood of, it may be advantageous to the discovery process to centre the parameter about this value. That is, one may introduce $\tilde{\mu} = \mu - \bar{\mu}$ and discover a mapping as a function of $\tilde{\mu}$ instead of $\mu$. The advantage this poses is that it can potentially maintain that the coefficients in the discovered mapping are kept relatively small and therefore the sparsity-promoting procedure can be more effective. To illustrate, consider the simple caricature  
	\begin{equation}
		(\mu - 10)^3 = \mu^3 + 30\mu^2 + 300\mu + 1000, 
	\end{equation}	 
where we can see that centring the parameter about $\bar{\mu} = 10$ results in a small single coefficient on $\tilde{\mu}^3 = (\mu - 10)^3$, whereas on the right the coefficients in powers of $\mu$ stretch across four orders of magnitude. Hence, not centring the parameter about $\bar{\mu} = 10$ has the potential to lead to numerical error and ill-conditioning of the problem. In what follows we will apply this centring procedure whenever applicable. 


\section{Full State Feedback Control of Orbits}\label{sec:LMI}

In this section, we will actively adjust the parameter value $\mu$ to control the behavior of the system using full-state feedback control~\cite{Dullerud,sp:book,Astrom2010book}. 
From the previous section, we may obtain a smooth, parsimonious, parameter-dependent mapping $\Fv:\R^d\times\R^p \to \R^d$ which faithful captures the iterates in a Poincar\'e section by 
\begin{equation}\label{Map}
	\xv_{n+1} = \Fv(\xv_n,\mu),
\end{equation}
for all $n \geq 0$. Assuming that at some $\bar{\mu} \in \R^p$ the map \eqref{Map} has a fixed point $\bar{\xv}\in\R^d$, i.e. $\bar{\xv} = \Fv(\bar{\xv},\bar{\mu})$, we may linearize the mapping about this fixed point to obtain 
\begin{equation}\label{LinearMap}
	\xv_{n+1} - \bar{\xv} \approx  \Av(\xv_n - \bar{\xv}) + \Bv(\mu_n - \bar{\mu})
\end{equation} 
where
\begin{equation}\label{ABMatrix}
	\Av := \Fv_x(\bar{\xv},\bar{\mu})\in\R^{d\times d},\quad  \Bv := \Fv_\mu(\bar{\xv},\bar{\mu})\in\R^{d\times p}.
\end{equation}
are the Jacobian matrices of \eqref{Map} resulting from differentiation with respect to $\xv$ and $\mu$, respectively. 
When $\xv_n$ is in a neighbourhood of $\bar{\xv}$, we will use the parameter fluctuations $\mu_n$ to control the system by making them linear functions of the state variable $\xv$, whereas when $\xv$ is far from $\bar{\xv}$ we will simply take $\mu_n = \bar{\mu}$. 
That is,  
\begin{equation}\label{paramControl}
	\mu_n = \begin{cases}
		\bar{\mu} + \Kv(\xv_n - \bar{\xv}) & |\xv_n - \bar{\xv}| \leq \eta \\
		\bar{\mu} &  |\xv_n - \bar{\xv}| > \eta 
	\end{cases} 
\end{equation} 
for a yet-to-be specified control matrix $\Kv\in\R^{p\times d}$ and some $\eta > 0$ which we refer to as the {\bf threshold parameter}. Notice that for $\xv_n$ in a neighbourhood of $\bar{\xv}$ with the parameter choices in \eqref{paramControl} the linearization \eqref{LinearMap} becomes
\begin{equation}\label{LinearMap2}
	\xv_{n+1}- \bar{\xv} \approx (\Av + \Bv\Kv)(\xv_n - \bar{\xv}),
\end{equation}	
and so our goal now becomes choosing the matrix $\Kv$ in such a way that every iterate of the linearized mapping \eqref{LinearMap2} converges to $\bar{\xv}$ as $n \to \infty$.  
This approach is known as \emph{full state feedback control} in the field of control theory, and there are well known solutions~\cite{Dullerud,sp:book,Astrom2010book,DDBook}. 

We will refer to a matrix $\Mv \in \R^{d\times d}$ as {\bf stable} if all eigenvalues $\lambda \in \mathbb{C}$ are such that $|\lambda| < 1$ and note that this implies that every solution of the linear mapping $\xv_{n+1} = \Mv \xv_n$ converges to zero as $n \to \infty$, regardless of initial condition. We make use of the following Lyapunov characterization which gives an equivalent condition for the stability of a matrix:
\begin{equation}\label{Stable}
	\Mv\mathrm{\ is\ stable} \iff \exists \Pv\succ 0\ \mathrm{s.t.}\ \Pv -  \Mv^T\Pv\Mv\succ 0. 
\end{equation}
Here the subscript $T$ denotes the matrix transpose and $\Pv \succ 0$ denotes that $\Pv$ is positive definite. Returning to \eqref{LinearMap2}, a matrix $\Kv$ making $(\Av + \Bv\Kv)$ stable is then equivalent to saying that there exists a $\Pv \succ 0$ such that
\begin{equation}\label{Stable2}
	\Pv - (\Av+\Bv\Kv)^T\Pv(\Av+\Bv\Kv) \succ 0.
\end{equation} 
Notice that to solve \eqref{Stable2} we are required to obtain both $\Kv$ and $\Pv$, thus making the left-hand-side nonlinear in the unknowns. In \S~\ref{subsec:FixedPt} we present Parrilo's method for obtaining a matrix $\Kv$ that makes $\Av+\Bv\Kv$ stable by converting \eqref{Stable2} to an equivalent LMI~\cite{Parrilo}. Such an LMI can be evaluated numerically as a semidefinite program. In \S~\ref{subsec:Cycle} we provide a simple extension of these methods to cyclic orbits.

Prior to presenting the equivalent LMI formulation for solving \eqref{Stable2}, we note that the process of obtaining a matrix $\Kv$ is referred to as the pole-placement method. Pole-placement is a standard approach in feedback control~\cite{Dullerud,Ogata,Sontag} to design a control matrix $\Kv$ to obtain arbitrary placement of the eigenvalues of the closed-loop system $\Av+\Bv\Kv$. 
Pole placement has already been proposed by Romeiras et al.~\cite{Romeiras} for the control of chaos. 
The arbitrary placement of the closed-loop eigenvalues is possible if and only if the system is \emph{controllable}~\cite{Dullerud,DDBook}, meaning that the Krylov subspace spanned by $\Av$ and $\Bv$ has full row rank.  
One common pole placement is obtained via a \emph{linear quadratic regulator} (LQR) that minimizes an infinite-time horizon cost function
\begin{equation}\label{CostFn}
	J = \sum_{n = 0}^\infty (\xv_n-\bar \xv)^T\Jv_1(\xv_n-\bar \xv) + (\mu_n-\bar\mu)^T\Jv_2(\mu_n-\bar\mu), \quad \Jv_1 \succeq 0, \Jv_2 \succ 0.
\end{equation}  
 In this case, the unique matrix $\Kv$ balances the aggressiveness of stabilization with control expenditure, and is obtained by solving a Riccati equation~\cite{DDBook}. 
 Notice that the value of $\eta$ can be used to control the value of any cost function in the form \eqref{CostFn} since the control of the nonlinear mapping \eqref{Map} is only turned on in a small neighbourhood of $\bar\xv$. Hence, we will not assume that a cost function is given since we may simply decrease the value of $\eta$ to decrease the strength of the control applied to the system. On the other end, the maximal allowable value of $\eta > 0$ is intimately related to the choice of $\Kv$ since it must be chosen to lie in the basin of attraction of $\bar{\xv}$ to the nonlinear mapping \eqref{Map} after implementing the parameter control \eqref{paramControl}. We therefore require that $|\Kv(\xv_n - \bar{\xv})|$ is sufficiently small to guarantee that $\mu_n$ remains in a neighbourhood of $\bar{\mu}$. We refer the reader to~\cite{Boccaletti,Romeiras} for a more complete discussion of the choice of the threshold parameter as it relates to the control matrix $\Kv$.


\subsection{LMIs for Stabilizing Fixed Points}\label{subsec:FixedPt}

Based on the preceding discussion, our goal is to obtain a matrix $\Kv$ such that $\Av+\Bv\Kv$ is stable using \eqref{Stable2}. Condition \eqref{Stable2} is nonlinear in the unknowns $\Pv$ and $\Kv$ and therefore we follow~\cite{Parrilo} to equivalently state this condition as an LMI for both $\Pv$ and $\Kv$. First, we use the Schur complement to write \eqref{Stable2} equivalently as 
\begin{equation}\label{StableControl2}
	\begin{bmatrix}
		\Pv & (\Av + \Bv\Kv)^T\Pv \\ \Pv(\Av+\Bv\Kv) & \Pv
	\end{bmatrix} \succ 0.
\end{equation}
Notice that condition \eqref{StableControl2} is only bilinear in $(\Kv,\Pv)$ as opposed to nonlinear in $\Kv$, but is still not linear in the unknown matrices. Since $\Pv \succ 0$, $\Pv$ is invertible, so we may define the symmetric matrix $\Qv = \Pv^{-1}$. Multiplying \eqref{StableControl2} on the left and right by the invertible block diagonal matrix ${\rm diag}(\Qv,\Qv)$ states \eqref{StableControl2} equivalently as 
\begin{equation}\label{StableControl3}
	\begin{bmatrix}
		\Qv & \Qv(\Av + \Bv\Kv)^T \\ (\Av+\Bv\Kv)\Qv & \Qv
	\end{bmatrix} \succ 0.
\end{equation}
Since \eqref{StableControl3} only contains terms of the form $\Kv\Qv$, we may define $\Yv = \Kv\Qv$ to arrive at the LMI in $(\Yv,\Qv)$ 
\begin{equation}\label{LMI}
	\begin{bmatrix}
		\Qv & \Qv\Av^T + \Yv^T\Bv^T \\ \Av\Qv + \Bv\Yv & \Qv
	\end{bmatrix} \succ 0,
\end{equation}
equivalent to solving \eqref{Stable2}. Notice that the Schur complement guarantees that if $(\Qv,\Yv)$ can be obtained to satisfy \eqref{LMI}, then $\Qv$ is necessarily invertible and so we recover $\Kv = \Yv\Qv^{-1}$. 

We now see that finding the control matrix $\Kv$ can be equivalently stated as an LMI and implemented as a semidefinite program. We may implement the strict positive definiteness of \eqref{LMI} numerically by considering some small $\varepsilon > 0$ and requiring that 
\begin{equation}\label{LMI_Numerical}
	\begin{bmatrix}
		\Qv - \varepsilon \Iv_d & \Qv\Av^T + \Yv^T\Bv^T \\ \Av\Qv + \Bv\Yv & \Qv - \varepsilon \Iv_d
	\end{bmatrix} \succeq 0,
\end{equation}
where $\Iv_d$ is the $d\times d$ identity matrix. In our applications in Sections~\ref{sec:Examples} and \ref{sec:ParamInd} we will take $\varepsilon = 10^{-3}$. Hence, we see that we have equivalently stated the controllability of the system \eqref{Map} through the parameter manipulations \eqref{paramControl} as a question of feasibility of a semidefinite program. That is, there exists a matrix $\Kv$ that makes $\Av+\Bv\Kv$ stable if and only if the convex cone of solutions $(\Qv,\Yv)$ to \eqref{LMI} is nonempty. One case for no such matrix $\Kv$ to exist is when $\Av$ has a left eigenvector $w$ such that $w^T\Av = \lambda w^T$, $w^T\Bv = 0$ and $|\lambda| \geq 1$. In such cases more sophisticated control techniques would be required to stabilize the fixed point $\bar{\xv}$ which are beyond the scope of this investigation.


\subsection{LMIs for Stabilizing Cycles} \label{subsec:Cycle}

We now extend the previous discussion of stabilizing fixed points of \eqref{Map} to stabilizing cycles. Let us assume that at $\mu = \bar{\mu}$ the set $\{\bar{\xv}_1,\bar{\xv}_2,\dots,\bar{\xv}_m\} \subset \R^d$, for some $m > 1$, satisfies 
\begin{equation}
	\bar{\xv}_{j+1} = \Fv(\bar{\xv}_j,\bar{\mu}), \quad j = 1,\dots,m
\end{equation}
where $\bar{\xv}_{m+1} = \bar{\xv}_1$. To stabilize the above $m$-cycle it is possible to move to the $m$th iterate map, $\Fv^m$, where each element of the $m$-cycle then becomes a fixed point and the analysis from the previous subsection applies. As pointed out in~\cite{Romeiras}, this would be suboptimal since the map $\Fv^m$ is overly sensitive to noise, especially in the present scenario when the mapping $\Fv$ is discovered from data and susceptible to numerical error. Therefore, we seek to apply parameter control similar to \eqref{paramControl} in a neighbourhood of each element of the $m$-cycle. 

We will again considering the iteration-dependent choice of the parameter, $\mu_n$, and so linearizing \eqref{Map} about each point on the $m$-cycle gives 
\begin{equation}\label{LinearMap3}
	\xv_{n+1} - \bar{\xv}_{j+1} \approx \Av_j(\xv_n - \bar{\xv}_j) + \Bv_j(\mu_n - \bar{\mu}),
\end{equation}   
where the matrices $\Av_j$ and $\Bv_j$ are defined as 
\begin{equation}
	\Av_j := \Fv_x(\bar{\xv}_j,\bar{\mu})\in\R^{d\times d},\quad  \Bv_j := \Fv_\mu(\bar{\xv}_j,\bar{\mu})\in\R^{d\times p}.
\end{equation}
for all $j = 1,\dots,m$. As before, we will select the parameter $\mu_n$ at each step of the iteration depending on the proximity of $\xv_n$ to an element of the $m$-cycle. We take 
\begin{equation}\label{paramControl2}
	\mu_n = \begin{cases}
		\bar{\mu} + \Kv_j(\xv_n - \bar{\xv}_j) & |\xv_n - \bar{\xv}_j| \leq \eta \\
		\bar{\mu} & \mathrm{otherwise} 
	\end{cases} 
\end{equation} 
for appropriately chosen matrices $\Kv_j\in\R^{p\times d}$. For $\xv_n$ in a neighbourhood of $\xv_j$ we have that \eqref{LinearMap3} becomes
\begin{equation}
	\xv_{n+1} - \bar{\xv}_{j+1} \approx (\Av_j + \Bv_j\Kv_j)(\xv_n - \bar{\xv}_j),	
\end{equation} 
for each $j = 1,\dots,m$. Linear stability of the $m$-cycle is then equivalent to obtaining matrices $(\Kv_1,\dots,\Kv_m)$ so that the matrix
\begin{equation}\label{CycleStable}
	(\Av_m + \Bv_m\Kv_m)\cdot(\Av_{m-1} + \Bv_{m-1}\Kv_{m-1})\cdots(\Av_2 + \Bv_2\Kv_2)\cdot(\Av_1 + \Bv_1\Kv_1)
\end{equation}   
is stable. We note that obtaining the $\Kv_j$ such that each $\Av_j + \Bv_j\Kv_j$ is stable is not sufficient to guarantee the matrix \eqref{CycleStable} is stable and so we require a method to handle the nonlinear Lyapunov condition \eqref{Stable} associated to the matrix \eqref{CycleStable}.  

To obtain sufficient matrices $(\Kv_1,\dots,\Kv_m)$ to make \eqref{CycleStable} stable we propose the following. Begin by obtaining a matrix $\Kv_m$ so that $\Cv_m:=\Av_m + \Bv_m\Kv_m$ is stable using the LMI \eqref{LMI}. Once an appropriate $K_m$ has been obtained, we obtain a matrix $\Kv_{m-1}$ so that the matrix
\begin{equation}
	\Cv_{m-1} := \Cv_m(\Av_{m-1} + \Bv_{m-1}\Kv_{m-1}) = \Cv_m\Av_{m-1} + \Cv_m\Bv_{m-1}\Kv_{m-1}
\end{equation}  
is stable. This can again be implemented as a semidefinite program using \eqref{LMI} with $\Cv_m\Av_{m-1}$ in the place of $\Av$ and $\Cv_m\Bv_{m-1}$ in the place of $\Bv$. We continue this process inductively by obtaining $\Kv_j$ so that the matrix 
\begin{equation}
	\Cv_j := \Cv_m\Cv_{m-1}\cdots \Cv_{j+1}(\Av_j + \Bv_j\Kv_j) = \Cv_m\Cv_{m-1}\cdots \Cv_{j+1}\Av_j + \Cv_m\Cv_{m-1}\cdots \Cv_{j+1}\Bv_j\Kv_j
\end{equation} 
is stable for all $j = 1,\dots,m-1$. Notice that the matrix $\Cv_1$ is exactly the matrix \eqref{CycleStable}, and therefore this process terminates in stabilizing the matrix \eqref{CycleStable} using the matrices $\Kv_j$. Moreover, at each step we may cast the problem of obtaining the $\Kv_j$ as an LMI which can be implemented numerically using \eqref{LMI_Numerical}, thus allowing for a principled automation of this process.


\section{Applications}\label{sec:Examples}

In this section we will provide illustrative examples of utility of our methods. All numerical procedures are implemented in MATLAB R2019a, while a robust python-based SINDy package (PySindy)~\cite{PySINDy} is publicly available at {\bf GitHub/dynamicslab/pysindy}. Throughout this section our library of candidate functions will consist of all monomials up to degree 5 in the Poincar\'e section variables. We use the MATLAB software YALMIP (version R20190425)~\cite{Yalmip} to solve the LMIs that determine the control matrices using the numerical implementation \eqref{LMI_Numerical}. All code for this section is available at {\bf GitHub/jbramburger/Stabilizing\_UPOs}. Finally, throughout this section we refer to period $m \geq 1$ orbits for continuous dynamical systems as periodic solutions which intersect the Poincar\'e section at $m$ distinct points, i.e. $m$-cycles. This is in contrast to the definition that the temporal period of the orbits is the integer $m$ which may be a source of confusion for the reader.


\subsection{The Chaotic H\'enon Map}

Our first example seeks to illustrate the stabilization process outlined in Section~\ref{sec:LMI} as it applies to a discrete dynamical system. We will focus on the H\'enon map
\begin{equation}\label{Henon} 
	\begin{split}
		x_{n+1} &= 1 - ax_n^2 + y_n, \\
		y_{n+1} &= bx_n,
	\end{split}
\end{equation}
in a neighbourhood of the standard chaotic parameter values $(a,b) = (1.4,0.3)$. The mapping \eqref{Henon} was originally presented as a simplified model of the Poincar\'e section of the Lorenz ordinary differential equations~\cite{Henon} and has grown into one of the canonical chaotic dynamical systems due to its rich dynamical structure. Most important to the discussion herein, it has been shown that the structure of the chaotic attractor can be understood in terms the unstable periodic orbits embedded in the attractor~\cite{HenonUPO} and therefore we will seek to stabilize these UPOs through slight adjustments to the parameter $b$ after each iteration. Throughout we will fix $a = 1.4$ for simplicity, and note the purpose of this section is to illustrate the effectiveness and the generalizability of the LMI method in \S~\ref{sec:LMI} for obtaining the control matrices, while not requiring one to potentially solve nonlinear equations as in~\cite{Wang}. 

\begin{figure} 
\center
\includegraphics[width = 0.49\textwidth]{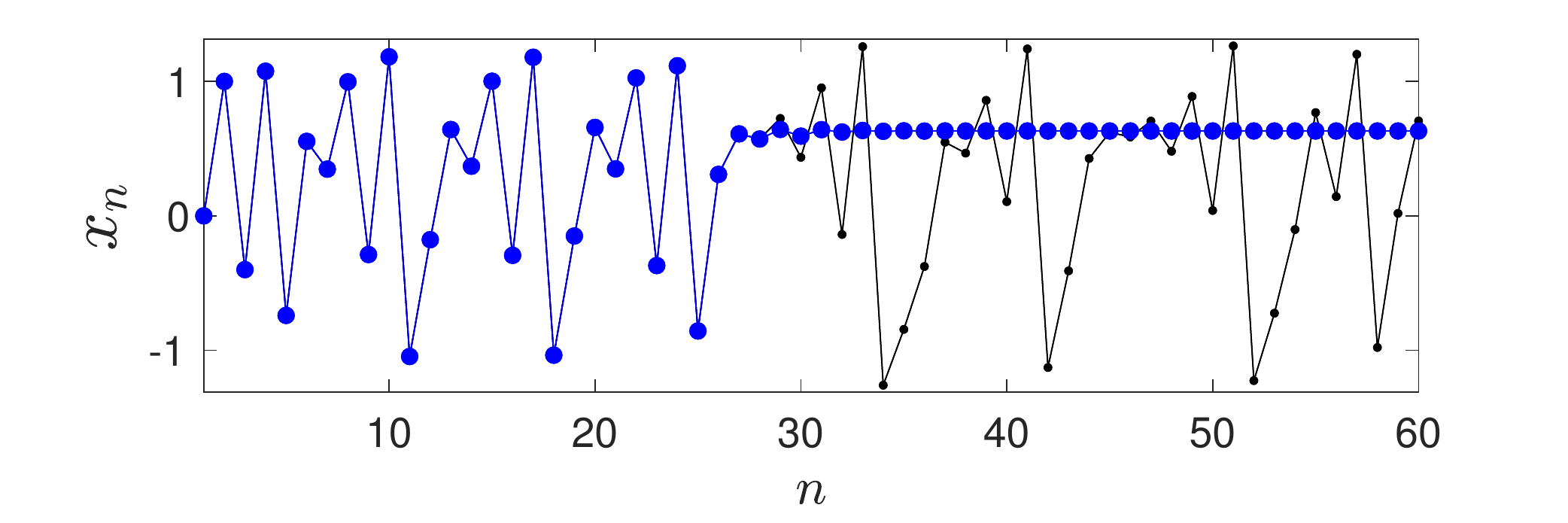} \ 
\includegraphics[width = 0.49\textwidth]{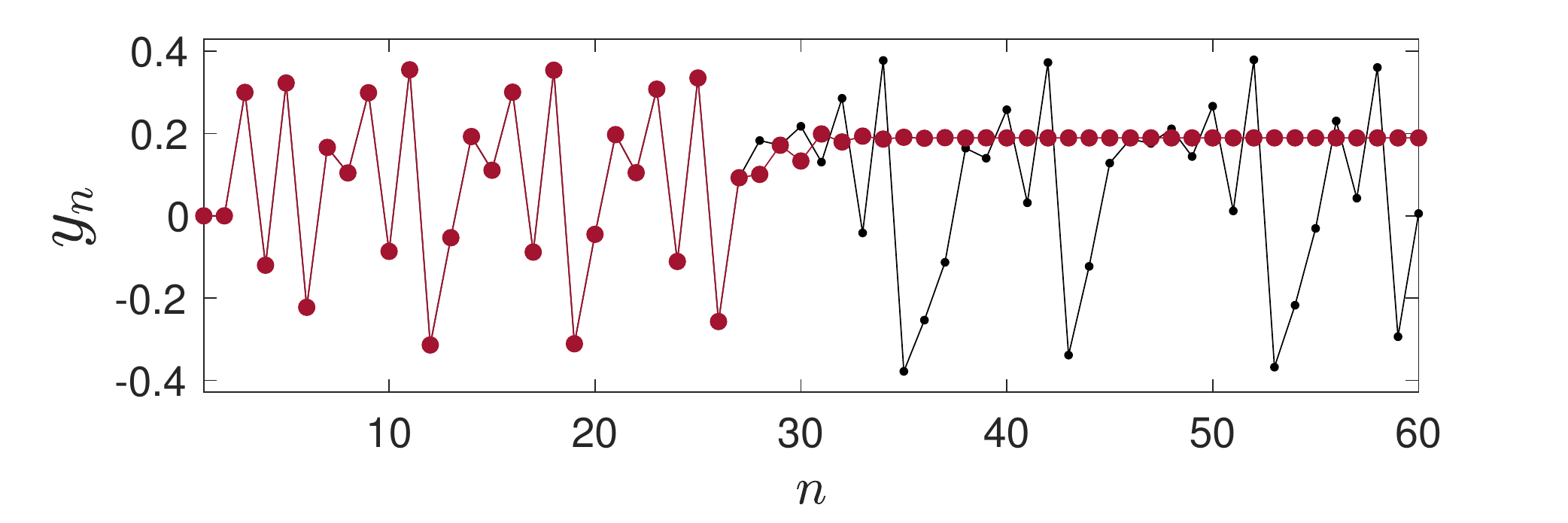}
\caption{Controlled versus uncontrolled orbits of the H\'enon map \eqref{Henon}. Controlled $x_n$ (left, blue) and $y_n$ (right, red) converge to the fixed point \eqref{HenonFixed}, whereas the uncontrolled orbits (black) wander chaotically along the attractor.}
\label{fig:Henon1}
\end{figure}

To begin, we note that at $(a,b) = (1.4,0.3)$ the fixed point
\begin{equation}\label{HenonFixed} 
	\begin{split}
		\bar{x} &= \frac{\sqrt{609} - 7}{28} \approx 0.63135 \\
		\bar{y} &= \frac{3(\sqrt{609} - 7)}{280} \approx 0.18941
	\end{split}
\end{equation}
of the mapping \eqref{Henon} lies in the chaotic attractor. Fixing $a = 1.4$ and considering $b$ in a neighbourhood of $0.3$ as our control parameter gives that the matrices $A \in \R^{2\times 2}$, $B\in\R^{2\times 1}$ from \eqref{ABMatrix} are given by
\begin{equation} 
	\Av = \begin{bmatrix}
		-2.8\bar{x} & 1 \\ 0.3 & 0	
	\end{bmatrix}, \quad 
	\Bv = \begin{bmatrix}
		0 \\ \bar{x}
	\end{bmatrix}.
\end{equation} 
We seek a matrix $\Kv \in \R^{1\times 2}$ so that $\Av + \Bv\Kv$ has eigenvalues entirely contained in the unit circle. Implementing \eqref{LMI} as a semidefinite program we obtain 
\begin{equation} 
	\Kv = \begin{bmatrix}
		-4.3796 & 2.3643
	\end{bmatrix}
\end{equation}	
which can be used to stabilize the fixed point \eqref{HenonFixed}. Indeed, the eigenvalues of $\Av + \Bv\Kv$ are approximately given by $-0.57639$ and $0.30129$, thus giving linearized stability of the fixed point. We present the resulting iterates of the controlled H\'enon mapping with a threshold parameter $\eta = 0.01$ and the uncontrolled mapping, both with initial conditions $(x_0,y_0) = (0,0)$, in Figure~\ref{fig:Henon1}. 

\begin{figure} 
\center
\includegraphics[width = 0.49\textwidth]{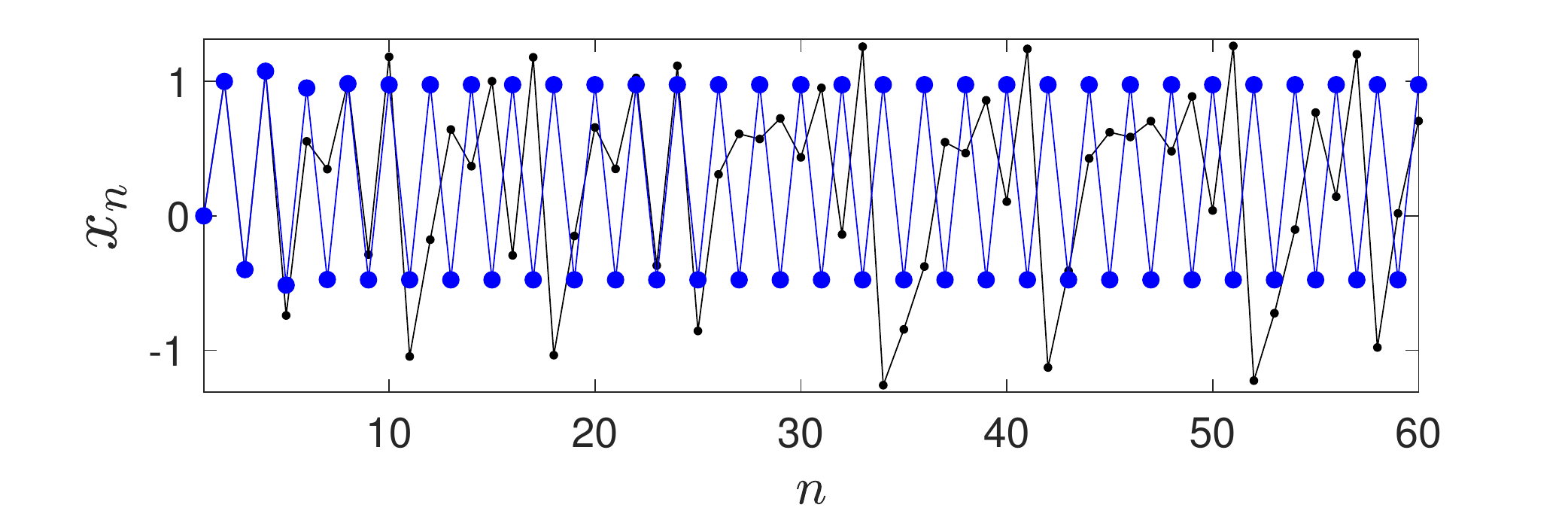} \ 
\includegraphics[width = 0.49\textwidth]{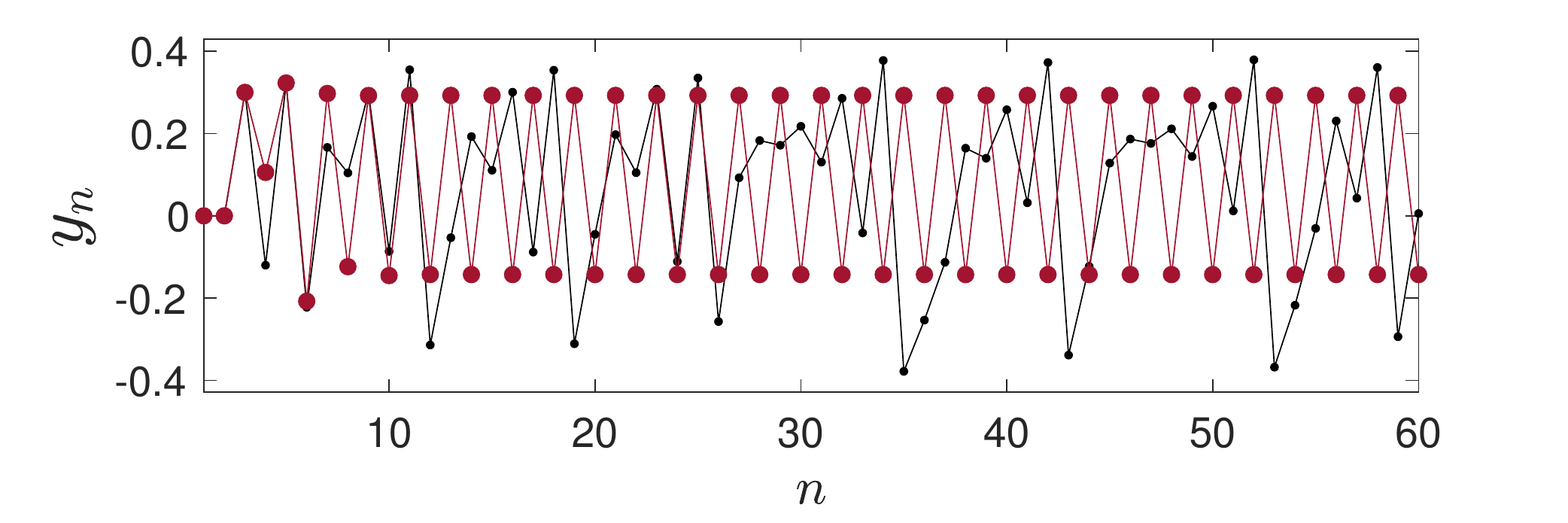}
\caption{Controlled versus uncontrolled orbits of the H\'enon map \eqref{Henon}. Controlled $x_n$ (left, blue) and $y_n$ (right, red) converge to a 2-cycle of \eqref{HenonFixed}, whereas the uncontrolled orbits (black) wander chaotically along the attractor.}
\label{fig:Henon2}
\end{figure}

We may similarly control the orbits of \eqref{Henon} to stabilize cyclic orbits of the map. For example, the H\'enon map has a 2-cycle given by
\begin{equation}\label{Henon2Cycle}
	\begin{split} 
		(x_1,y_1) &= (-0.47580,0.29274)  \\
		(x_2,y_2) &= (0.97580,-0.14274)
	\end{split}
\end{equation}  
for which the resulting linearized matrices are given by
\begin{equation} 
	\Av_{1,2} = \begin{bmatrix}
		-2.8\bar{x}_{1,2} & 1 \\ 0.3 & 0		
	\end{bmatrix} \quad 
	\Bv_{1,2} = \begin{bmatrix}
		0 \\ \bar{x}_{1,2}
	\end{bmatrix}.
\end{equation}
The procedure of \S~\ref{subsec:Cycle} then calls for obtaining matrices $\Kv_{1,2}\in\R^{1\times 2}$ such that the eigenvalues of $\Av_2 + \Bv_2\Kv_2$ and $(\Av_2 + \Bv_2\Kv_2)(\Av_1 + \Bv_1\Kv_1)$ are contained within the unit circle of the complex plane. Our numerical implementation provides the matrices 
\begin{equation} 
	\Kv_1 = \begin{bmatrix} 
		-7.1378 & 2.5333
	\end{bmatrix}, \quad \Kv_2 = \begin{bmatrix} 
		-6.8982  -5.6484
	\end{bmatrix},
\end{equation}
which give that the eigenvalues of $\Av_2 + \Bv_2\Kv_2$ are $-0.45567$ and $0.19546$, while the eigenvalues of $(\Av_2 + \Bv_2\Kv_2)(\Av_1 + \Bv_1\Kv_1)$ are $-0.07729$ and $-0.00204$. The controlled orbits are presented in Figure~\ref{fig:Henon2}, where we again use a threshold parameter of $\eta = 0.01$.  

\begin{figure} 
\center
\includegraphics[width = 0.8\textwidth]{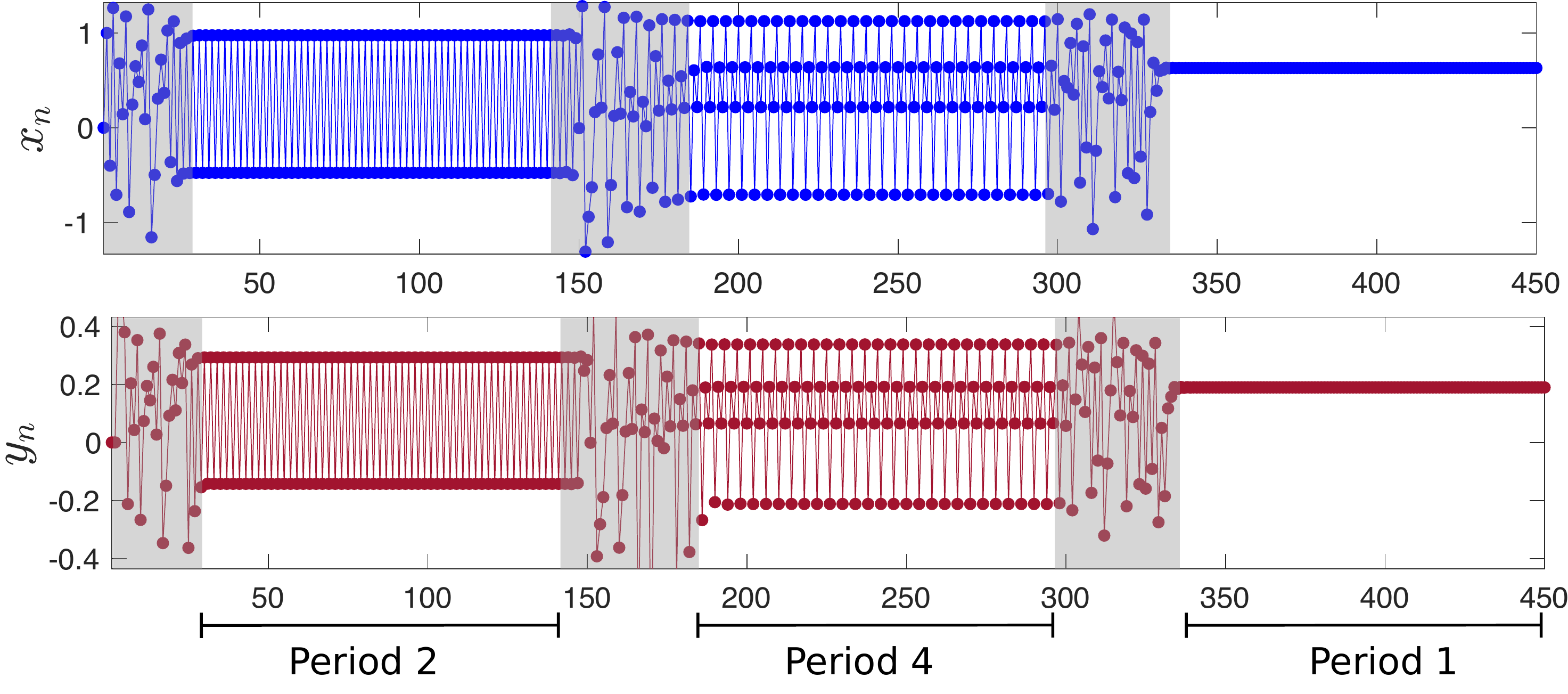} 
\caption{Controlled switching orbits $x_n$ (top, blue) and $y_n$ (bottom, red) of \eqref{HenonFixed}. Orbits are forced to a 2-cycle for 100 iterates, then a 4-cycle for 100 iterates, and then to a fixed point. Shaded areas represent transient orbits between the forced cyclic motion where control is not active.}
\label{fig:HenonSwitch}
\end{figure}

This process may be continued to stabilize cyclic orbits of the H\'enon map \eqref{Henon} of any period, as well as allow for switching between periodic orbits. We present in Figure~\ref{fig:HenonSwitch} a controlled orbit that switches from a 2-cycle to a 4-cycle to a fixed point. In this case we have taken the threshold parameter to be $\eta = 0.05$ for the simple reason that for smaller $\eta$ the orbit with initial conditions $x(0) = y(0) = 0$ experiences a longer transient that makes the visual presentation cluttered. The orbits are forced to remain near each of the periodic orbits for 100 iterates.

\subsection{Isolated Periodic Orbits}\label{subsec:Hopf}

Here will we now introduce the data-driven discovery of Poincar\'e maps aspect of our method to a simple toy model that exhibits an isolated UPO. The system of interest is given in cartesian coordinates by
\begin{equation}\label{Hopf} 
	\begin{split}
		\dot{x} &= -\omega y + x(x^2 + y^2 - \mu^2)(16 - x^2 - y^2) \\
		\dot{y} &= \omega x + y(x^2 + y^2 - \mu^2)(16 - x^2 - y^2)
	\end{split}
\end{equation}	
where $\omega > 0$ and $\mu \in (0,4)$ are system parameters. System~\eqref{Hopf} is more conveniently written in polar coordinates $(r,\theta)$ with $x = r\cos(\theta)$ and $y = r\sin(\theta)$, resulting in the system
\begin{equation}\label{HopfPolar} 
	\begin{split}
		\dot{r} &= r(r^2 - \mu^2)(16 - r^2) \\
		\dot{\theta} &= \omega. 
	\end{split}
\end{equation}
From the polar form \eqref{HopfPolar} it is easy to see that for each $\mu \in (0,4)$ the origin of system \eqref{Hopf} is stable and has two limit cycles given by the circles centred at the origin with radius $\mu$ and $4$. The larger limit cycle is stable for all $\mu \in (0,4)$, whereas the smaller limit cycle is unstable in the same parameter range. In what follows we will centre our analysis about $\mu = 2$ and show that our methods can be applied to keep trajectories of \eqref{Hopf} close to the circle $x^2 + y^2 = 4$ through slight variations in the parameter $\mu$ in a neighbourhood of $2$.  

We will take our Poincar\'e section to be the half-line $\{(x,y)|\ y = 0,\ x\geq0\}$, which can equivalently be stated using the polar equation \eqref{HopfPolar} by the line $\theta = 0$. Notice that since the radial and azimuthal components of system \eqref{HopfPolar} decouple, there is a fixed time of $2\pi/\omega$ that trajectories take to return to the section. When $\omega > 0$ is small this return time is very long and therefore intersections of trajectories with the Poincar\'e section will exhibit large jumps when near the unstable orbit $x^2 + y^2 = \mu^2$. To circumvent this, we will fix $\omega = 100\pi$, thus allowing us to gather sufficient section training data near the unstable limit cycle by slightly slowing the divergence from it. Then, to generate the section data we sweep $\mu$ from $1$ to $3$ in increments of $0.1$ and for each $\mu$ we generate five trajectories with the following initial conditions: $x(0) = 0,\mu/2,\mu,(\mu + 4)/2,4$ and $y(0) = 0$ in all cases. Writing $x_n$ to be the iterates in the Poincar\'e section and applying the mapping discovery method of \S~\ref{sec:SINDy} results in the mapping
\begin{equation}\label{HopfPSec} 
	\begin{split}
		x_{n+1} = 0.34953&x_n -0.31661x^2_n  -0.32565(\mu - 2)x_n  + 0.52203x_n^3  -0.32058(\mu - 2)x_n^2 \\ 
		& -0.10056x_n^4  + 0.10050(\mu - 2)x_n^3.
	\end{split}
\end{equation}   
At $\mu = 2$ the mapping \eqref{HopfPSec} has three nonnegative fixed points given to four decimal places by $x = 0,2,4$, corresponding to the closed orbits of the system \eqref{Hopf}. Using the notation \eqref{ABMatrix} we have $\Av = 2.1295$ and $\Bv = -1.1296$, and taking $\Kv \in(0.99991,2.7704)$ guarantees that $\Av + \Bv\Kv\in(-1,1)$, thus stabilizing the fixed point $x = 2$ of \eqref{HopfPSec} at the parameter value $\mu = 2$. 

\begin{figure} 
\center
\includegraphics[width = 0.49\textwidth]{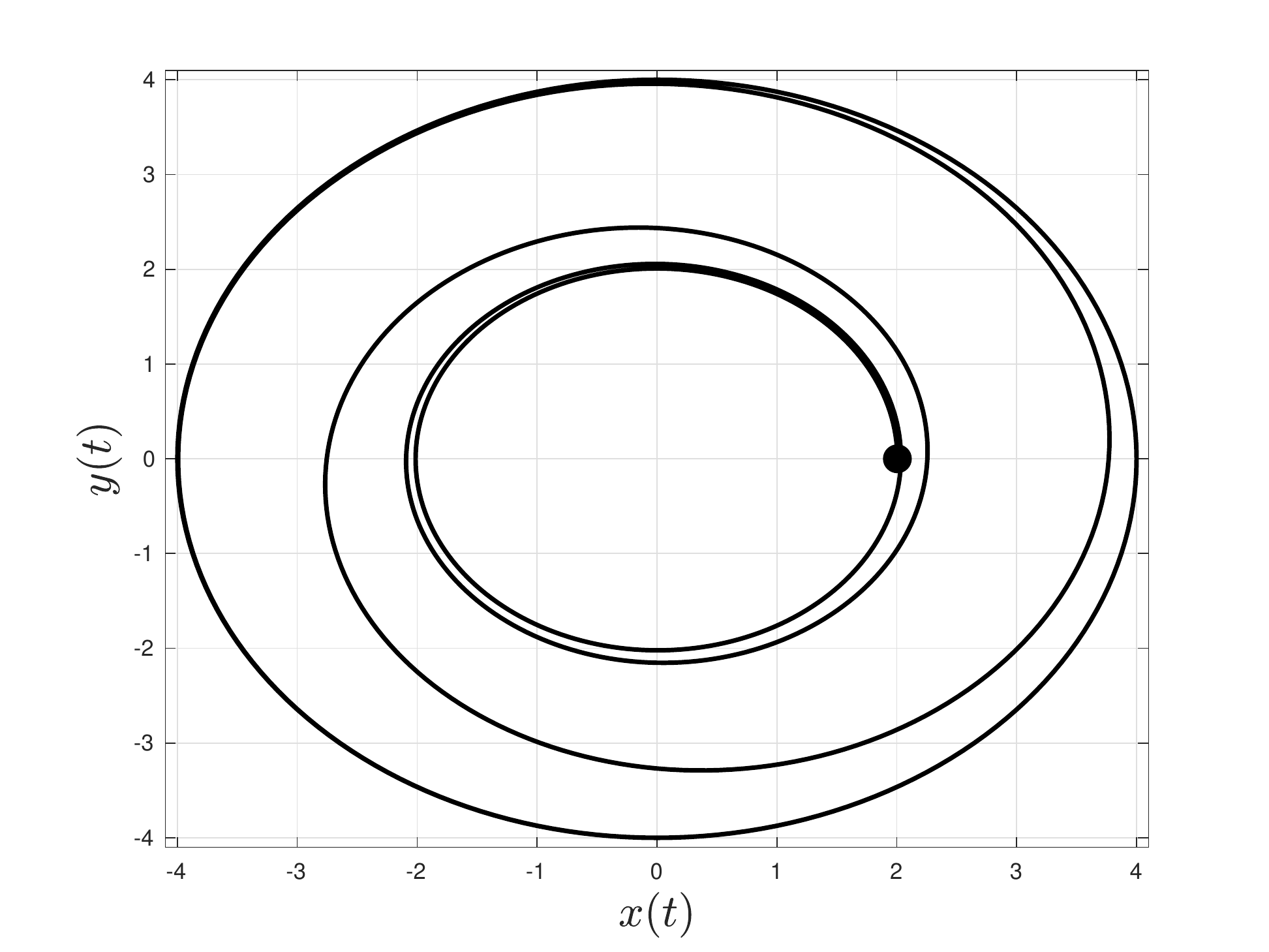} \ 
\includegraphics[width = 0.49\textwidth]{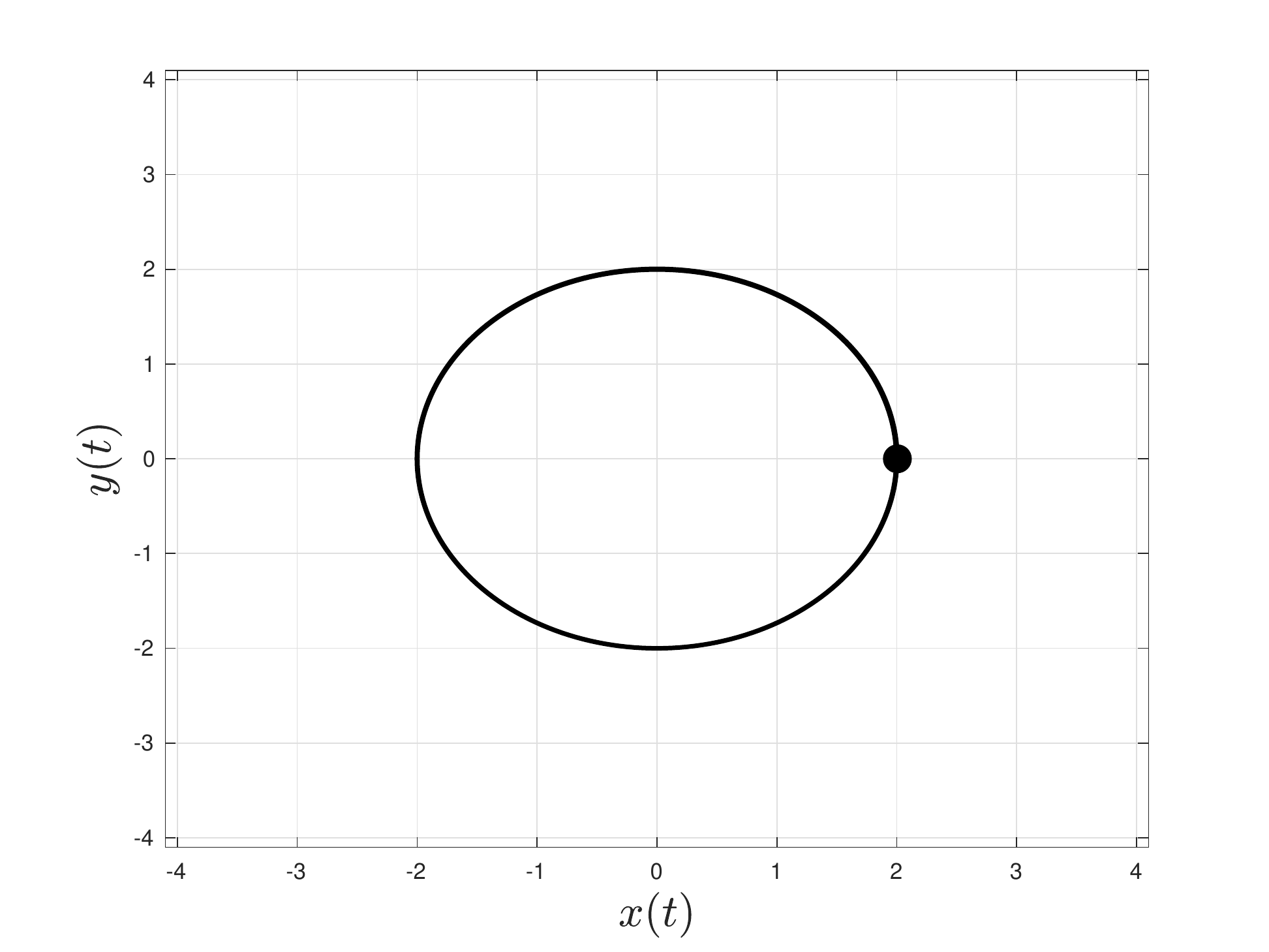}
\caption{Uncontrolled (left) and controlled (right) trajectories of the system \eqref{Hopf}. Both trajectories have initial condition $(x(0),y(0)) = (2.005,0)$, indicated by a large black dot. On the left the uncontrolled orbit slowly diverges from the UPO $x^2 + y^2 = 4$ and converges to the stable orbit $x^2 + y^2 = 16$. On the right we take $\Kv = 1.3$ and $\eta = 0.01$ to guarantee that the trajectory remains close to the UPO for all time.}
\label{fig:Hopf}
\end{figure}

In Figure~\ref{fig:Hopf} we present the results of our stabilization procedure to give a trajectory of \eqref{Hopf} that remains close to the UPO $x^2 + y^2 = 4$ with $\mu$ in a neighbourhood of $2$. We also present the results of numerical integration of the system without control for reference. Our stabilization procedure uses $\Kv = 1.3$, a threshold parameter $\eta = 0.01$, and the controlled trajectory provided in Figure~\ref{fig:Hopf} has initial condition $(x(0),y(0)) = (2.005,0)$. These values are all intimately related since the above analysis requires $|\Kv(x_n - 2)|$ is sufficiently small and bounded by $\eta$. Clearly increasing $\Kv$ requires that $|x_n - 2|$ decreases for our linear stability analysis to remain valid, and so increasing $\Kv$ will require $|x(0) - 2|$ to decrease.

\subsection{The R\"ossler System}\label{subsec:Rossler}

\begin{figure} 
\center
\includegraphics[width = 0.49\textwidth]{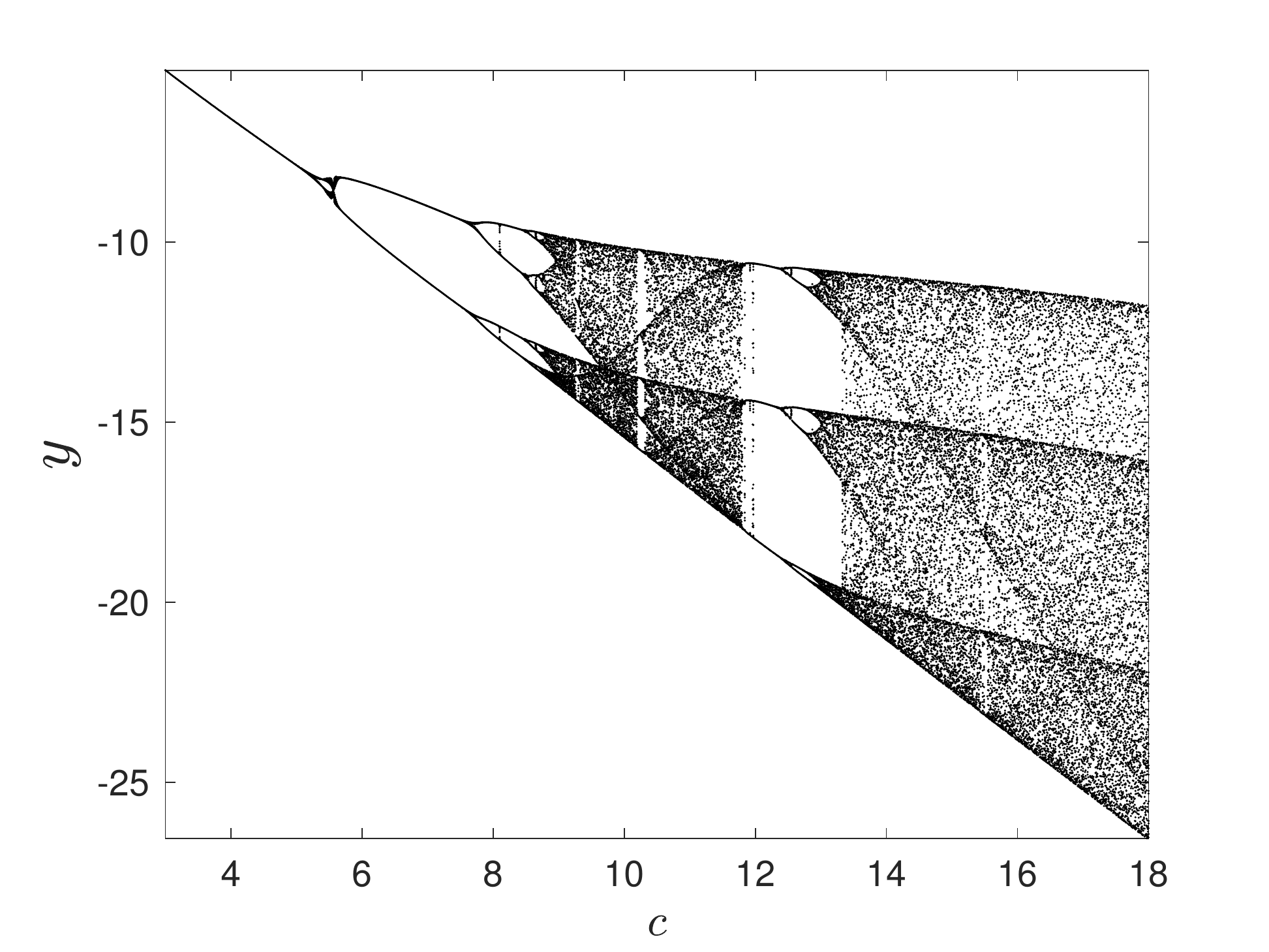} 
\caption{Poincar\'e section data for the R\"ossler system \eqref{Rossler} for varying $c$.}
\label{fig:Rossler_BifDiag}
\end{figure}

In this subsection we will now apply our methods to stabilizing UPOs of the three-dimensional R\"ossler system, given by~\cite{Rossler} 
\begin{equation}\label{Rossler} 
	\begin{split}
		\dot{x} &= -y - z \\
		\dot{y} &= x + 0.1y \\
		\dot{z} &= 0.1 + z(x-c)
	\end{split}
\end{equation}
where $c \in \R$ is our control parameter. We refer the reader to the works~\cite{Barrio,Peitgen} for detailed discussions of the bifurcations and structure of the attractor in the system. To summarize, the R\"ossler system is well-known for its sequence of period-doubling bifurcations of periodic orbits that lead to a chaotic attractor as $c$ is increased. Moreover, the system has the property that trajectories always cross the $x = 0$ plane at $z = 0$. Therefore, we will take our Poincar\'e section to be when trajectories cross $x = 0$ from negative to positive, thus allowing us to consider parameter-dependent Poincar\'e maps that depend only on the value of $y$ in these sections. In Figure~\ref{fig:Rossler_BifDiag} we plot data from this Poincar\'e section for $c \in [3,18]$ for reference throughout this section.

\begin{figure} 
\center
\includegraphics[width = 0.49\textwidth]{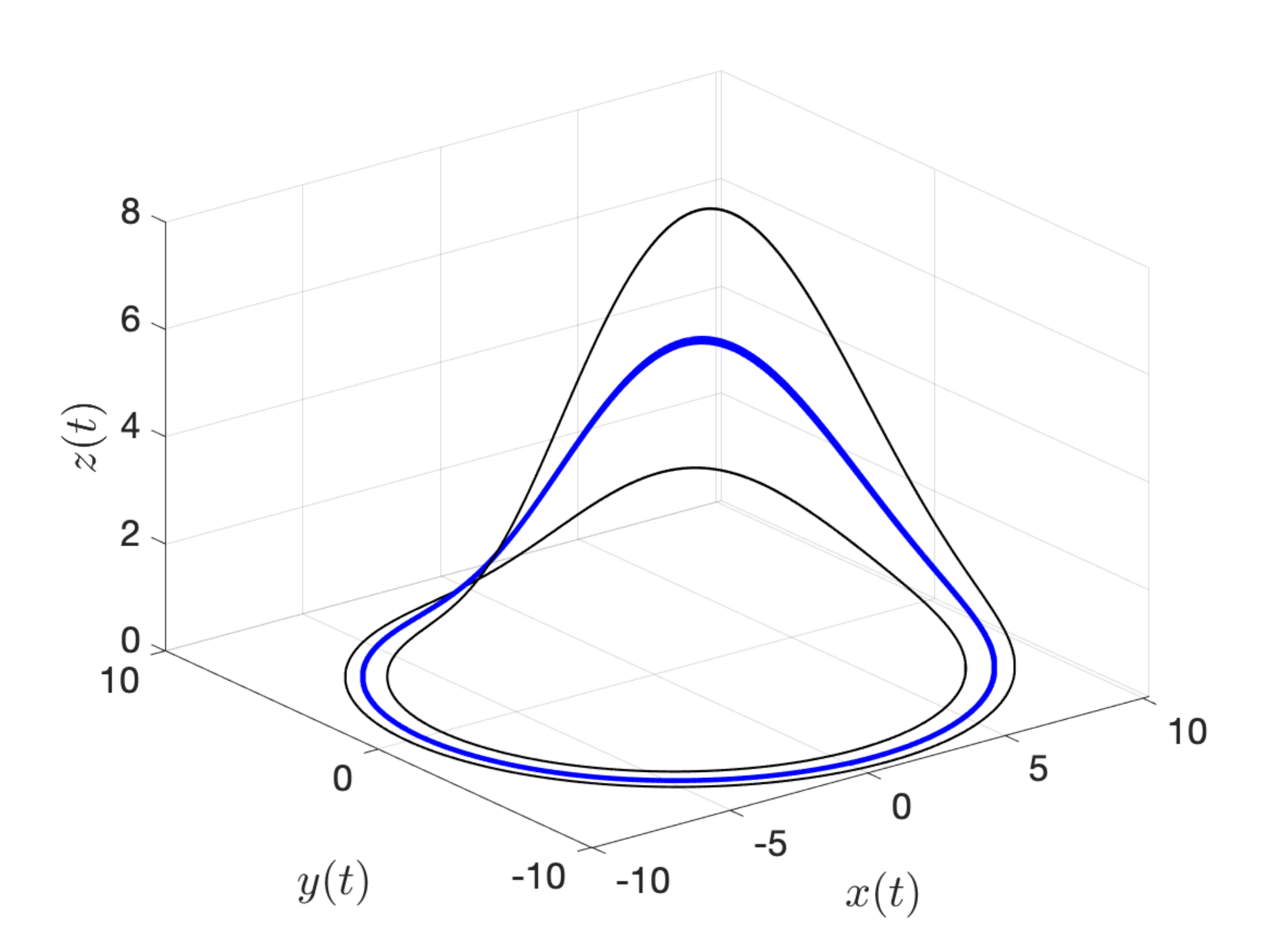} \ 
\includegraphics[width = 0.49\textwidth]{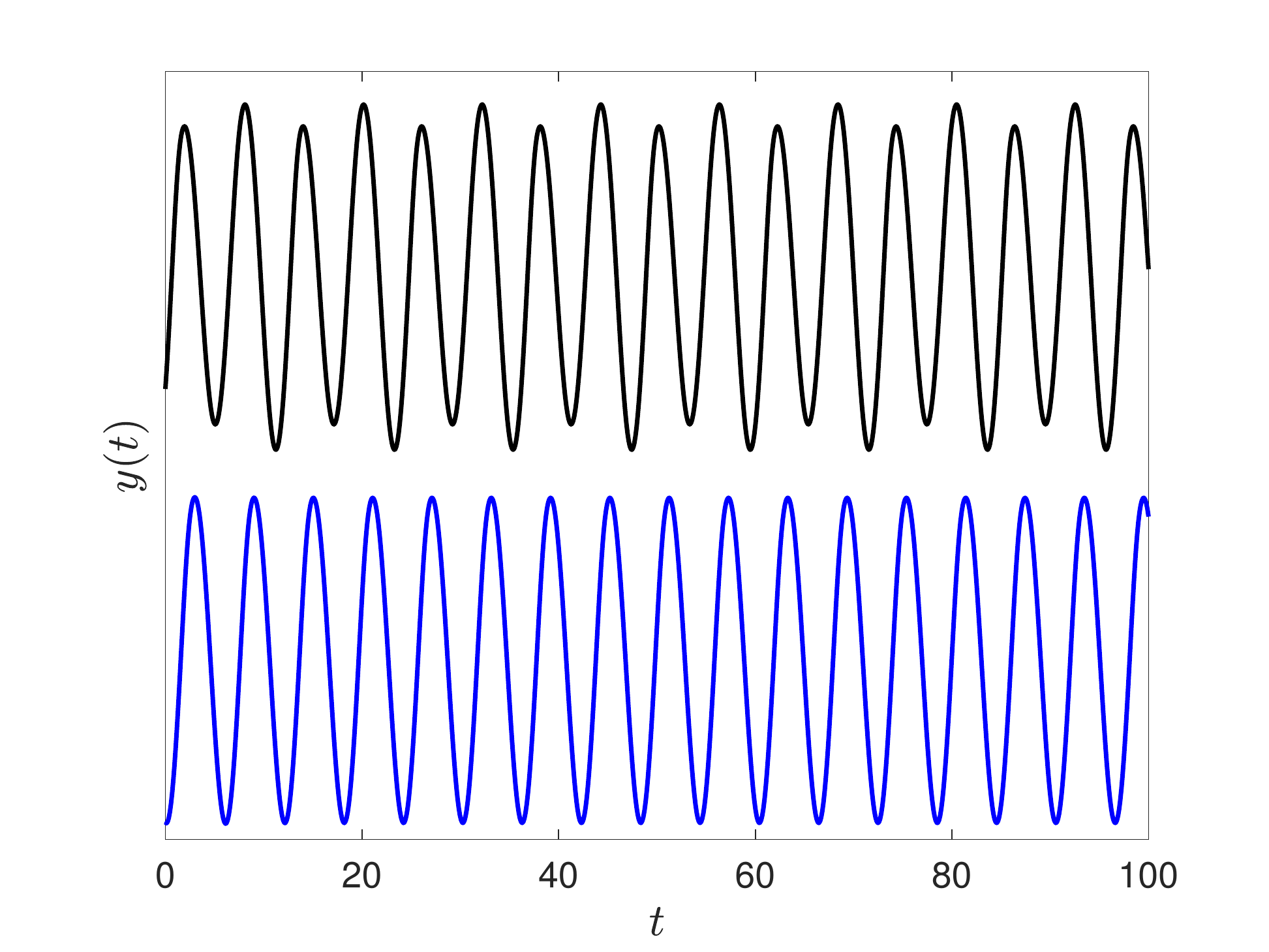}
\caption{Left: A stabilized period 1 orbit (blue) of the R\"ossler system \eqref{Rossler} with $c = 6$. Also plotted is the period 2 attractor (black). Right: The $y(t)$ components of the two periodic orbits to emphasize the periods of oscillation.}
\label{fig:Rossler6}
\end{figure}

We begin by centring our analysis about $c = 6$, where from Figure~\ref{fig:Rossler_BifDiag} we see that the attractor is a period 2 orbit. We can further see that this period 2 orbit is birthed from a period-doubling bifurcation when a period 1 orbit destabilizes and therefore at $c = 6$ we expect there to be a period 1 UPO. To find and stabilize this period 1 UPO at $c = 6$ we gather section data with $c$ progressing from 3 to 8 in increments of $0.1$ by initializing a single trajectory of \eqref{Rossler} with initial condition $(x(0),y(0),z(0)) = (0,-5,0)$ and integrating forward to $t = 100$ to provide sufficient section data for each $c$. In this parameter range our training data sees the period doubling bifurcation take place and therefore we expect the discovered mapping to be able to track the period 1 orbit as it destabilizes in the period-doubling bifurcation. Indeed, using a sparsity parameter of $\lambda = 0.01$ the discovered mapping is given by
\begin{equation}\label{c6}
	\begin{split}
	y_{n+1} &= 1.0607 + 2.1337y_n  + 0.84685c  + 0.41696y_n^2  + 0.73883cy_n  + 0.26510c^2 -0.044678y_n^3 \\  &-0.16800cy_n^2  -0.20357c^2y_n  -0.076355c^3
	\end{split}
\end{equation}	 
which has the unstable fixed point $y = -9.1238$ at $c = 6$. From the mapping \eqref{c6} we find that we require $\Kv\in\R$ so that
\begin{equation}
	-1.1343 -2.6562\Kv \in (-1,1)
\end{equation} 
to stabilize the period 1 UPO. In Figure~\ref{fig:Rossler6} we present the results of our stabilization procedure with $(\Kv ,\eta) = (0.5,0.1)$ on a trajectory with initial condition starting in the Poincar\'e section with $y(0) = -9.1338$. We also provide the period 2 attractor for comparison. 

\begin{figure} 
\center
\includegraphics[width = 0.49\textwidth]{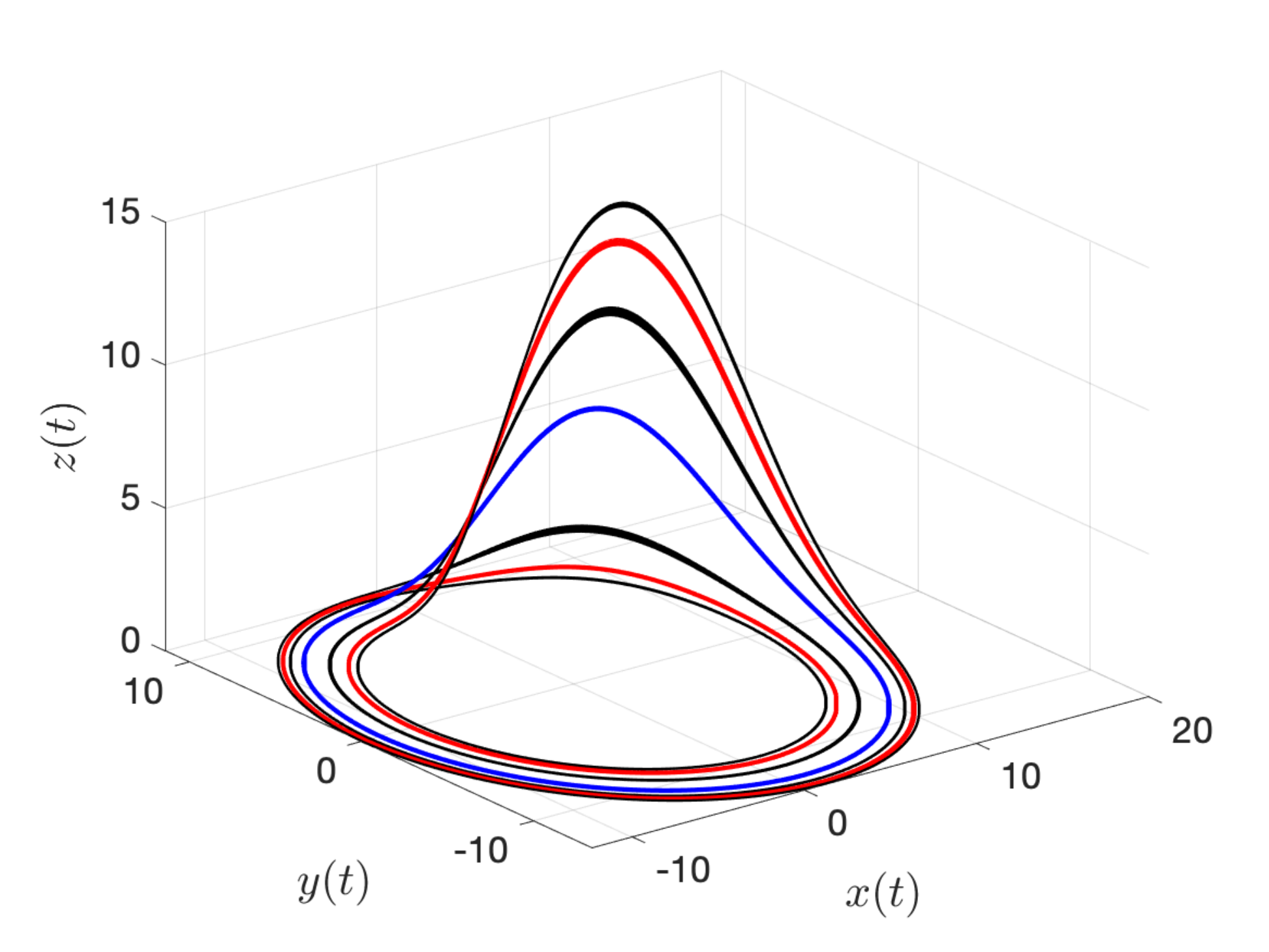} \ 
\includegraphics[width = 0.49\textwidth]{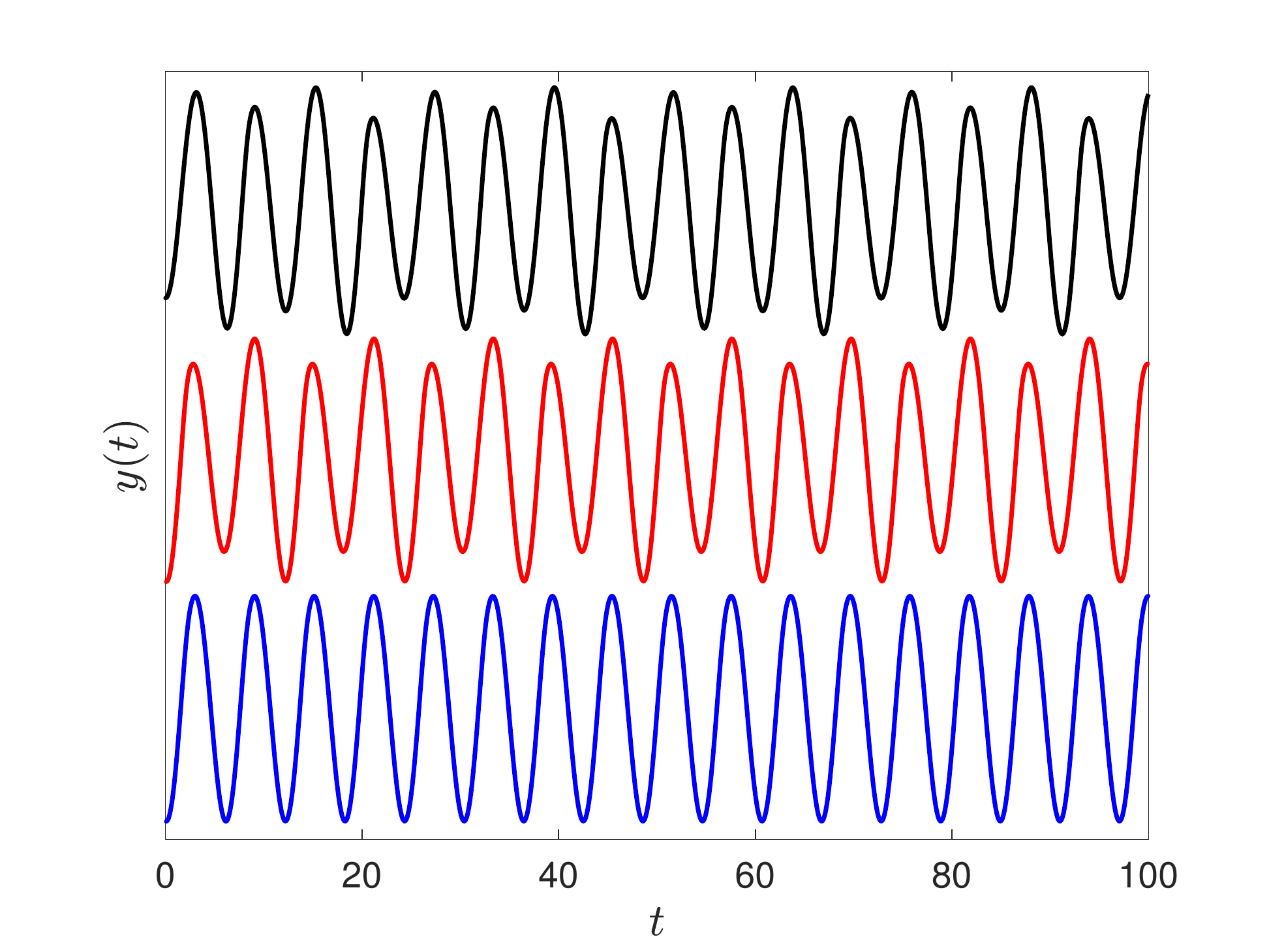}
\caption{Left: A stabilized period 1 orbit (blue) and period 2 orbit (red) of the R\"ossler system \eqref{Rossler} with $c = 8.5$. Also plotted is the period 4 attractor (black). Right: The $y(t)$ components of the three periodic orbits to emphasize the periods of oscillation.}
\label{fig:Rossler85}
\end{figure}

It is apparent from Figure~\ref{fig:Rossler_BifDiag} that for some $c \in (7,8)$ the period 2 orbit undergoes another period doubling bifurcation, thus destabilizing and resulting in a stable period 4 orbit for $c$ slightly above this bifurcation point. From this we expect that there is both a period 1 and a period 2 UPO at $c = 8.5$, which we show can be stabilized using our methods. Here we gather Poincar\'e section data for $c$ running from 6 to 9, noting that only the second period-doubling bifurcation is present in the training data. Hence, this training data does not explicitly include any information on the period 1 orbit. Using a sparsity parameter of $\lambda = 0.01$ we find the resulting mapping is given by
\begin{equation}\label{c8.5}
	\begin{split}
	y_{n+1} &= 5.16011 + 3.7081y_n  + 1.408c  -0.19798y_n^2  -1.5633cy_n  + -1.5431c^2  + -0.21689y_n^3 -0.68975cy_n^2 \\ &-0.57559c^2y_n -0.054612c^3  + -0.010789y_n^4 -0.032203cy_n^3  + -0.024004c^2y_n^2,
	\end{split}
\end{equation} 
which for $c = 8.5$ has an unstable fixed point at $y = -12.135$ and an unstable 2-cycle given by the sequence $-13.104 \rightarrow -10.122 \rightarrow -13.104$. Using the mapping we can again find appropriate control values $\Kv \in \R$ to stabilize the period 1 and period 2 orbits present at $c=8.5$. A resulting numerical integration of this stabilization procedure is presented in Figure~\ref{fig:Rossler85}, with the blue trajectory representing the period 1 orbit and the red trajectory having period 2. We also plot the stable period 4 orbit for reference. 

\begin{figure} 
\center
\includegraphics[width = 0.49\textwidth]{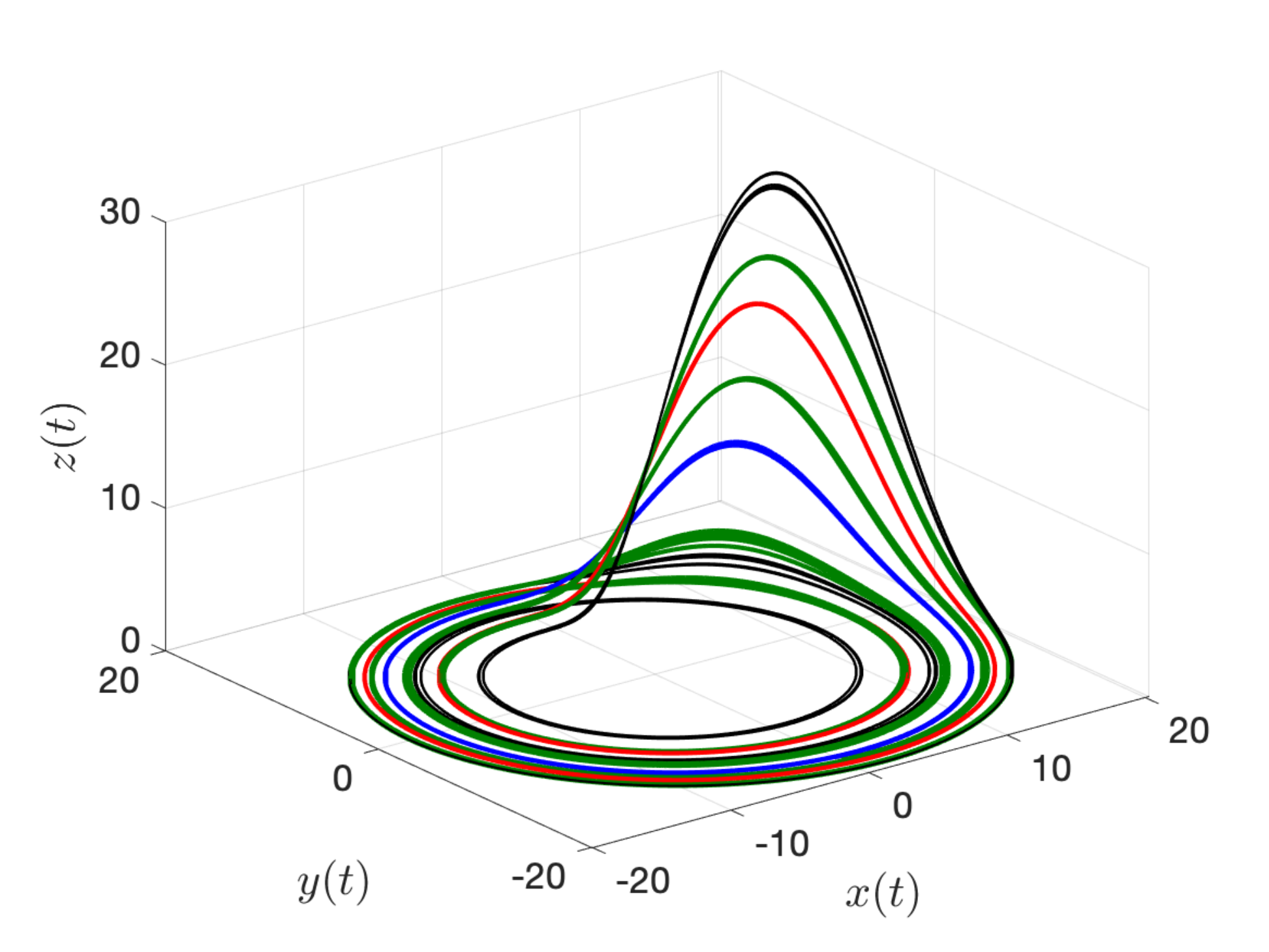} \ 
\includegraphics[width = 0.49\textwidth]{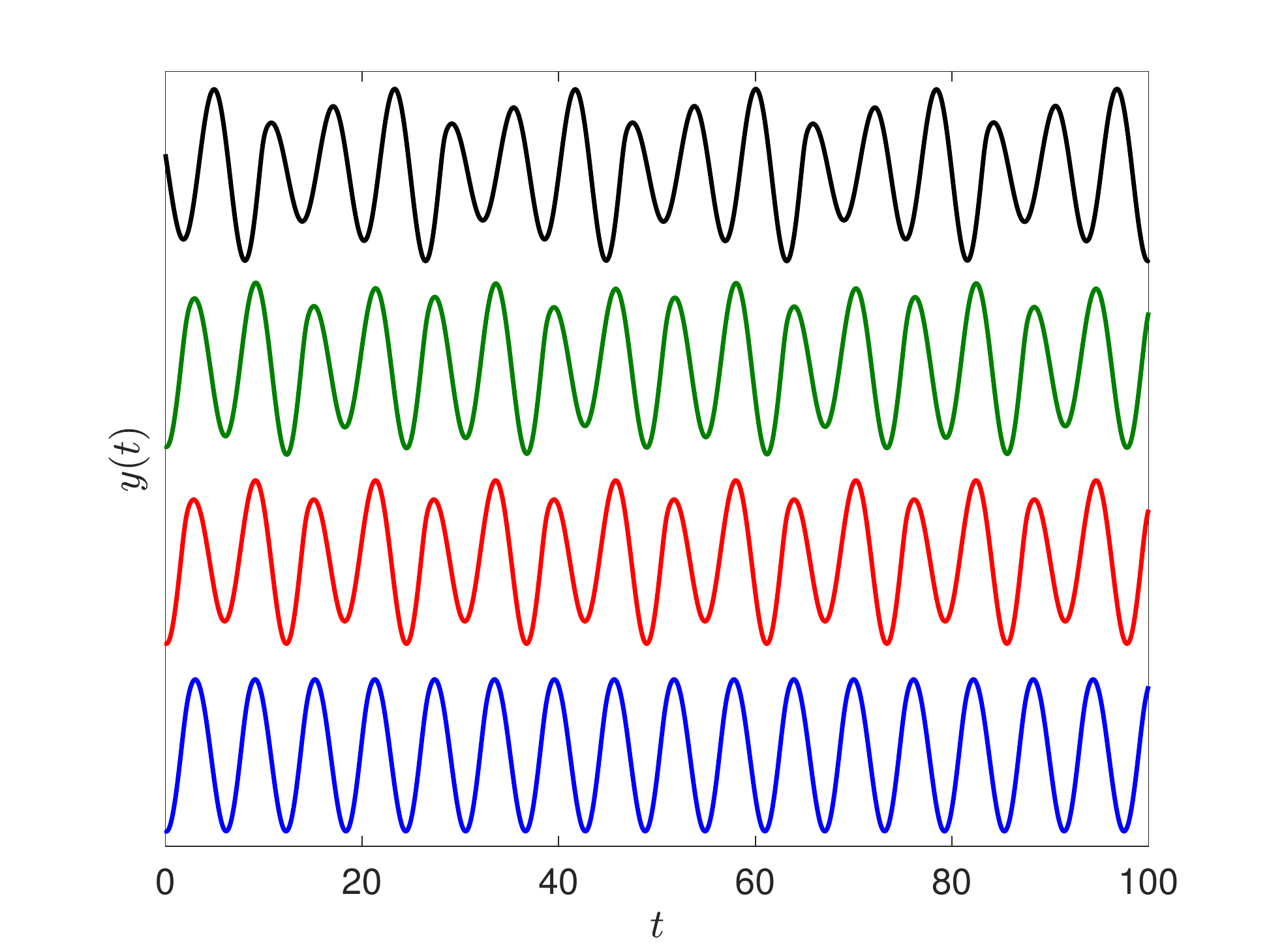}
\caption{Left: Stabilized period 1 (blue), period 2 (red), and period 4 (green) orbits  of the R\"ossler system \eqref{Rossler} with $c = 12.6$. Also plotted is the period 6 attractor (black). Right: The $y(t)$ components of the four periodic orbits to emphasize the periods of oscillation.}
\label{fig:Rossler126}
\end{figure}

So far we have seen how we may use destabilizing bifurcations to build information into the discovered Poincar\'e mapping. We now turn to a more challenging example where the training data provided does not include information about the destabilizing bifurcations. To this end we focus on the parameter value $c = 12.6$, which from Figure~\ref{fig:Rossler_BifDiag} we see that there is a period 6 attractor, but we also expect there to be a chaotic invariant set since we are beyond the initial period-doubling cascade into chaos for the system. The method of mapping discovery is analogous to the previous examples where we drag the parameter $c$ from $12.4$ to $12.8$ in increments of $0.01$ and use a sparsity parameter $\lambda = 0.01$, resulting in the mapping 
\begin{equation}\label{c12.6}
	\begin{split}
		y_{n+1} &= -24.619 -6.3815y_n -6.1923(c - 12.6)  -0.78238y_n^2  -1.3190(c - 12.6)y_n -11.943(c - 12.6)^2 \\  &-0.025552y_n^3 -0.066161(c - 12.6)y_n^2 -2.2118(c - 12.6)^2y_n  -74.692(c - 12.6)^3 - 0.081830(c - 12.6)^2y_n^2 \\ &-10.927(c - 12.6)^3y_n  + 198.72(c - 12.6)^4 -0.46384(c - 12.6)^3y_n^2  + 25.007(c - 12.6)^4y_n \\ & + 325.66(c - 12.6)^5,
	\end{split}
\end{equation}  
where we have centred the parameter dependence about $c = 12.6$ to minimize coefficients. Importantly, iterations of the mapping \eqref{c12.6} exhibits the following unstable periodic sequences
\begin{equation}
	\begin{gathered}
		-16.897 \rightarrow -16.897 \\
		-18.125 \rightarrow -13.569 \rightarrow -18.124 \\
		-10.859 \rightarrow -14.782 \rightarrow -19.041 \rightarrow -10.859 \\
		-12.758 \rightarrow -17.489 \rightarrow -15.632 \rightarrow -18.441 \rightarrow -12.758 
	\end{gathered}	
\end{equation} 
which constitute a fixed point, 2-cycle, 3-cycle, and 4-cycle, respectively. In Figure~\ref{fig:Rossler126} we plot the resulting period 1,2, and 4 orbits using the mapping \eqref{c12.6} and our stabilization procedure, along with the stable period 6 orbit. Not shown is the stabilized period 3 orbit since it is nearly indistinguishable from the period 6 orbit when plotted in $(x,y,z)$-space, but we note that by comparing the numerical values along the period 3 and period 6 trajectories our method presented herein does in fact stabilize this orbit as well. Furthermore, it may be possible to use \eqref{c12.6} to find higher order periodic orbits of the R\"ossler system \eqref{Rossler}, although numerical error becomes compounded through successive compositions of \eqref{c12.6} and therefore the resulting values of the higher order cycles may become inaccurate as the order of the cycle increases.   

As a final example for the R\"ossler system, we turn our attention to the parameter value $c = 18$, where from Figure~\ref{fig:Rossler_BifDiag} we expect the attractor to be chaotic. It was shown in~\cite{Bramburger} that in the case of a fully chaotic attractor the discovered mappings should not be expected to replicate the training data since slight numerical error can have drastic effects on the dynamics. Our goal here is to stabilize a period 1 and period 2 orbit present for $c = 18$ using training data from $c = 17$ to $c = 19$ in increments of $0.01$. We centre our parameter at $c = 18$ and use a sparsity parameter of $\lambda = 10^{-4}$ to obtain the mapping
\begin{equation}\label{c18}
	\begin{split}
		y_{n+1} &= 87.357 + 35.111y_n  + 18.527(c - 18)  + 4.8787y_n^2  + 5.9000(c - 18)y_n  + 0.53593(c - 18)^2 \\ &+ 0.33036y_n^3  + 0.64487(c - 18)y_n^2  + 0.14988(c - 18)^2y_n  + 0.067907(c - 18)^3  + 0.010451y_n^4 \\  &+ 0.028726(c - 18)y_n^3  + 0.013416(c - 18)^2y_n^2  -0.027(c - 18)^3y_n  + 0.20912(c - 18)^4 \\ & + 0.0001225y_n^5  + 0.0004345(c - 18)x^4  + 0.00027937(c - 18)^2y_n^3  -0.00034414(c - 18)^3y_n^2 \\ & + 0.03308(c - 18)^4y_n  -0.27141(c - 18)^5.
	\end{split}
\end{equation} 
We comment that using such a small sparsity parameter has the effect that many terms are included in the mapping that have very small coefficients, but we find that these terms are necessary to accurately obtain the location of the fixed point and 2-cycle in the Poincar\'e section. Precisely, using a sparsity parameter of $\lambda \geq 10^{-3}$ results in fixed points and 2-cycles that differ by $\mathcal{O}(1)$ to those found in \eqref{c18}, having the effect that we cannot stabilize the periodic orbits using these values due to their numerical inaccuracy.

\begin{figure} 
\center
\includegraphics[width = 0.49\textwidth]{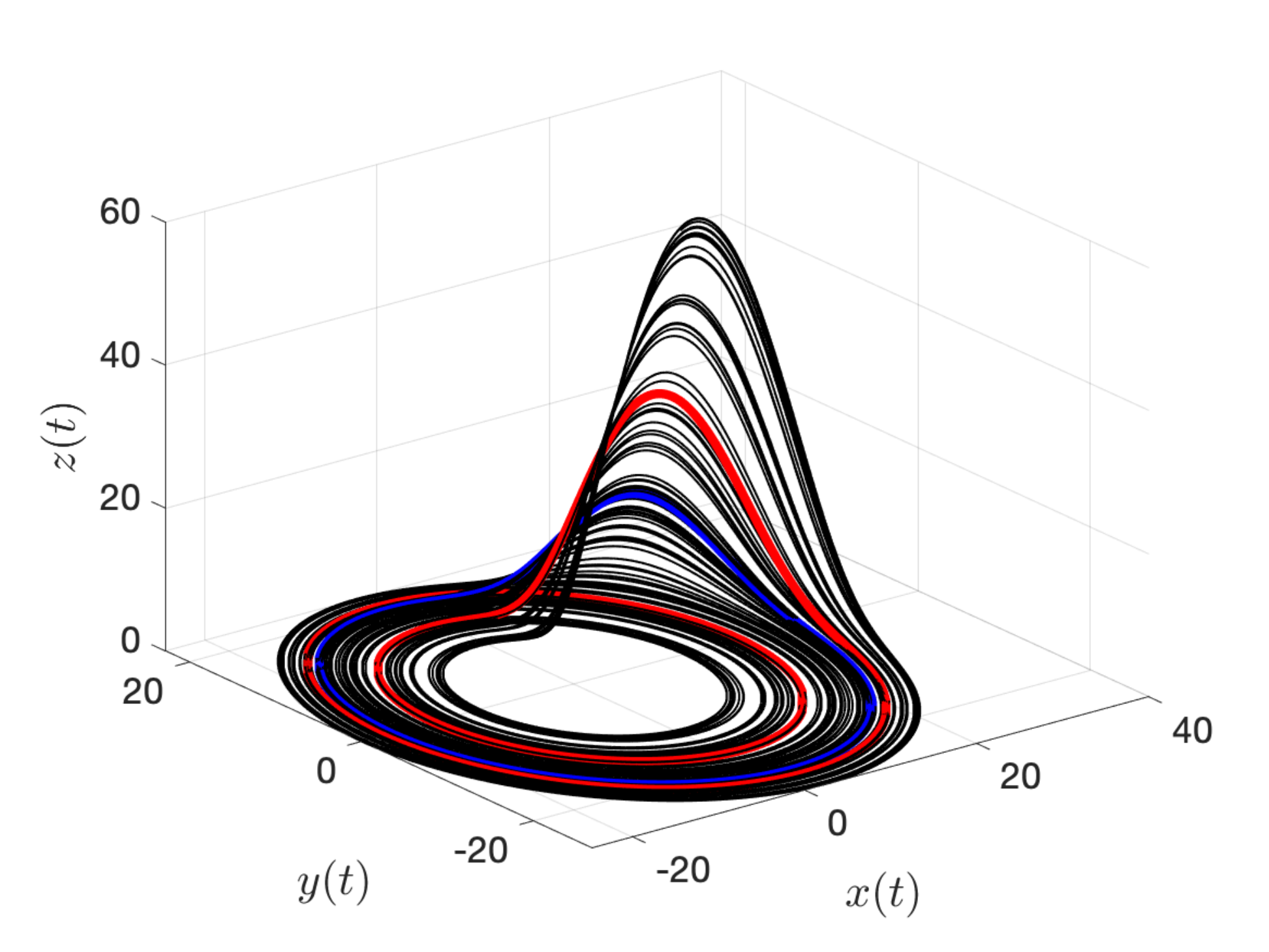} \ 
\includegraphics[width = 0.49\textwidth]{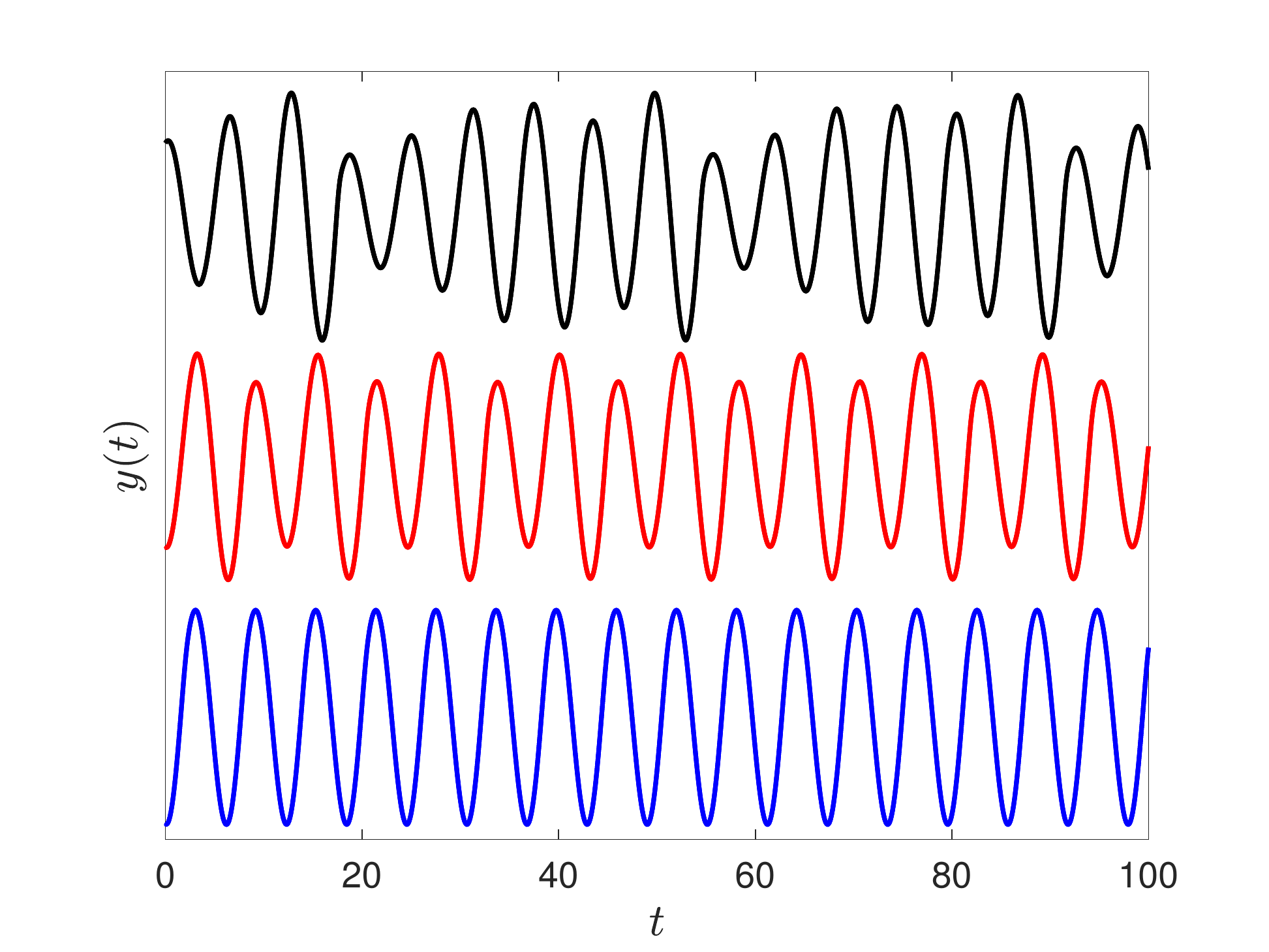}
\caption{Left: Stabilized period 1 (blue) and period 2 (red) orbits of the R\"ossler system \eqref{Rossler} along with the chaotic attractor (black) at $c = 18$. Right: The $y(t)$ components of the orbits to emphasize the periods of oscillation.}
\label{fig:Rossler18}
\end{figure}

At $c = 18$ the mapping \eqref{c18} has an unstable fixed point $y = -22.905$ and unstable 2-cycle $-17.802 \rightarrow -24.749 \rightarrow -17.802$. In Figure~\ref{fig:Rossler18} we present the stabilized period 1 and 2 orbits along with the uncontrolled chaotic attractor. Control parameters are given by $(\Kv,\eta) = (-0.6,0.1)$ for the period 1 orbit and $(\Kv_1,\Kv_2,\eta) = (0.5,-0.55,1)$ for the period 2 orbit. We note that the period 2 orbit requires a larger threshold parameter $\eta$ to stabilize the orbit. This could be attributed to numerical error in the values of the 2-cycle presented above. Furthermore, we report that we were unable to stabilize a period 3 orbit, potentially due again to numerical error compounded through successive compositions of the mapping \eqref{c18}, leading to inaccurate values of the 3-cycle. We hope to overcome this issue in a follow-up study.        

\subsection{Sprott's Chaotic Jerk System}\label{subsec:Sprott}

In this example we will apply our method to a system whose Poincar\'e section is two-dimensional. Consider the chaotic jerk system of Sprott given by the third-order equation
\begin{equation}\label{Sprott1}
	\frac{\drm^3 x}{\drm t^3} + \mu\frac{\drm^2 x}{\drm t^2} - \bigg(\frac{\drm x}{\drm t}\bigg)^2 + x = 0. 
\end{equation}
Equation \eqref{Sprott1} can equivalently be written as the first-order system of equations
\begin{equation}\label{Sprott2}
	\begin{split}
		\dot{x} &= y, \\
		\dot{y} &= z, \\
		\dot{z} & = - x -\mu z + y^2,
	\end{split}
\end{equation}
where the dot represents differentiation with respect to $t$. Much like the R\"ossler differential equation, system \eqref{Sprott2} experiences a sequence of period-doubling bifurcations leading to chaos as one decreases the parameter $\mu$ from $2.15$ down to $2$. We refer the reader to~\cite{Sprott} for a more complete discussion of the system \eqref{Sprott2} and particularly to Figure~1 of~\cite{Sprott} for the full bifurcation diagram. We simply note that the bifurcation diagram of \eqref{Sprott2} bears a significant resemblance to that of the R\"ossler system presented in Figure~\ref{fig:Rossler_BifDiag}. 

Our Poincar\'e section will be when trajectories cross $y = 0$ from positive to negative. Unlike the R\"ossler system, the dynamics in this section are genuinely two-dimensional, thus adding a level of complexity to the stabilization of UPOs in \eqref{Sprott2}. For brevity we will focus exclusively on stabilizing period 1 orbits near various choices of the bifurcation parameter $\mu$. The parameter choices are $\mu = 2.1, 2.06,$ and $2.05$ where the attractor is a period 2 orbit, a period 8 orbit, and chaotic, respectively. Our results are summarized in Figure~\ref{fig:Sprott} with the stabilized period 1 orbit presented in blue along with the attractor in black. The Poincar\'e mapping near $\mu = 2.1$ was discovered from data generated with $\mu$ running from $2.08$ to $2.12$ in increments of $0.001$ and training data at each parameter generated from a single trajectory with $(x(0),y(0),z(0)) = (-6,0,2.5)$. We use the same Poincar\'e mapping to find and stabilize the period 1 orbits at $\mu = 2.05,2.06$ by gathering training data again from a single trajectory with the same initial conditions as before and letting $\mu$ run from $2.04$ to $2.08$ in increments of $0.001$.     

\begin{table} 
\center
\begin{tabular}{>{\centering\arraybackslash}m{0.32\textwidth} >{\centering\arraybackslash}m{0.32\textwidth}  >{\centering\arraybackslash}m{4cm}}
\includegraphics[width = 0.35\textwidth]{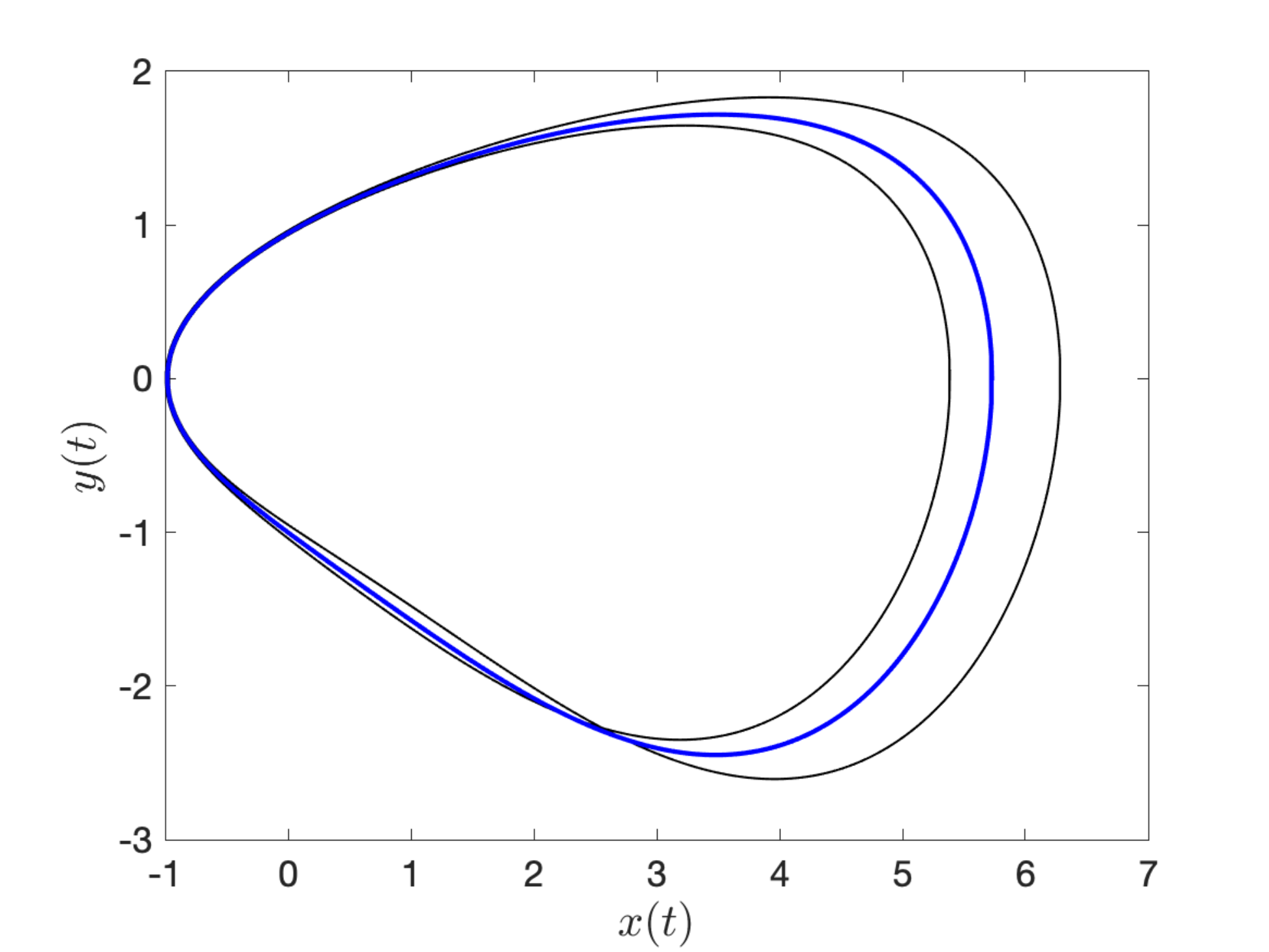} & 
\includegraphics[width = 0.35\textwidth]{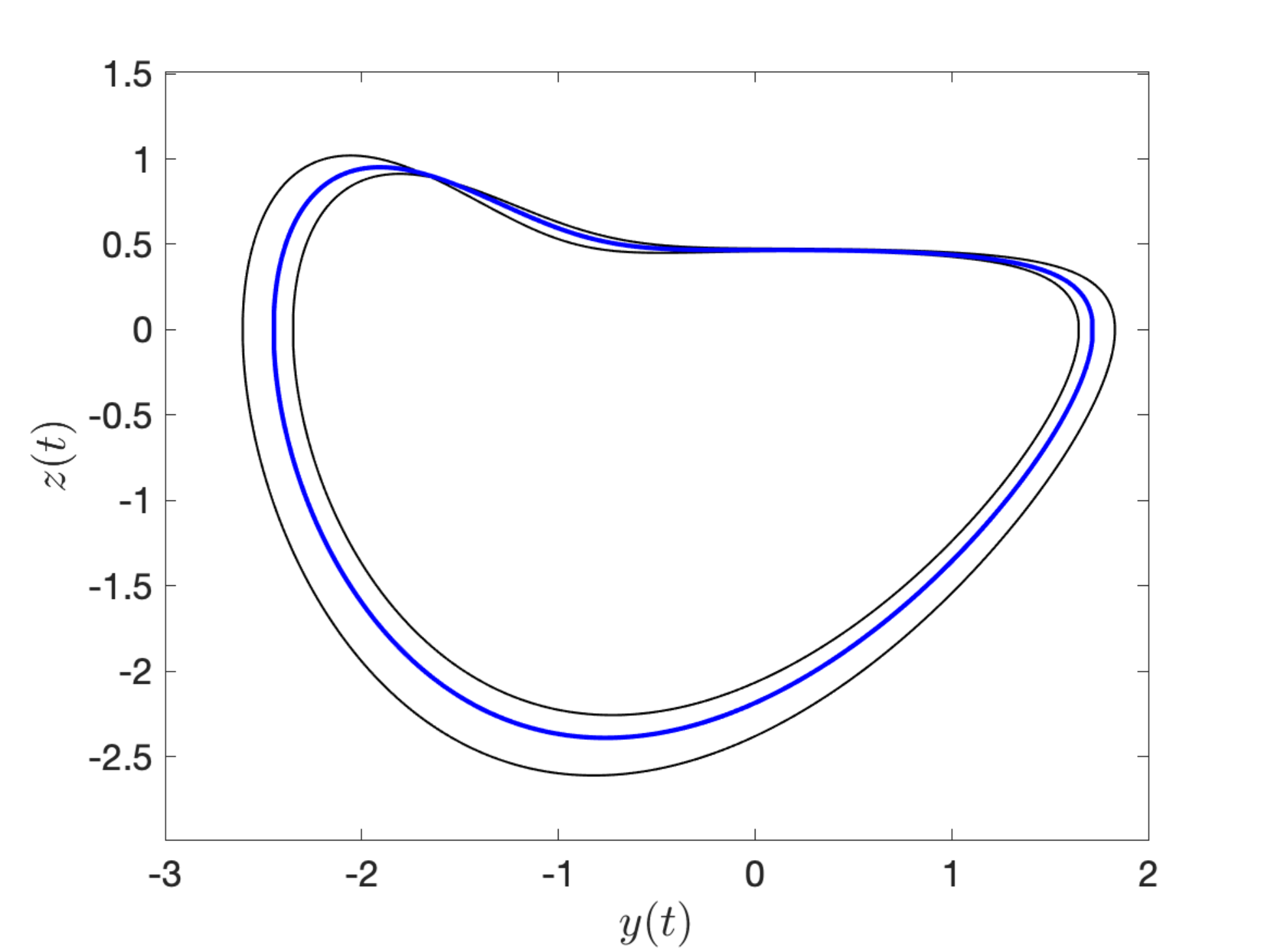} & \begin{tabular}{c} $\bar{\mu} = 2.1$ \\ $(\bar{x},\bar{z}) = (5.7043,-2.1278)$ \\ $\Kv = [0.11082\ 0.0077518]$ \\ $\eta = 0.1$ \end{tabular} \\
\includegraphics[width = 0.35\textwidth]{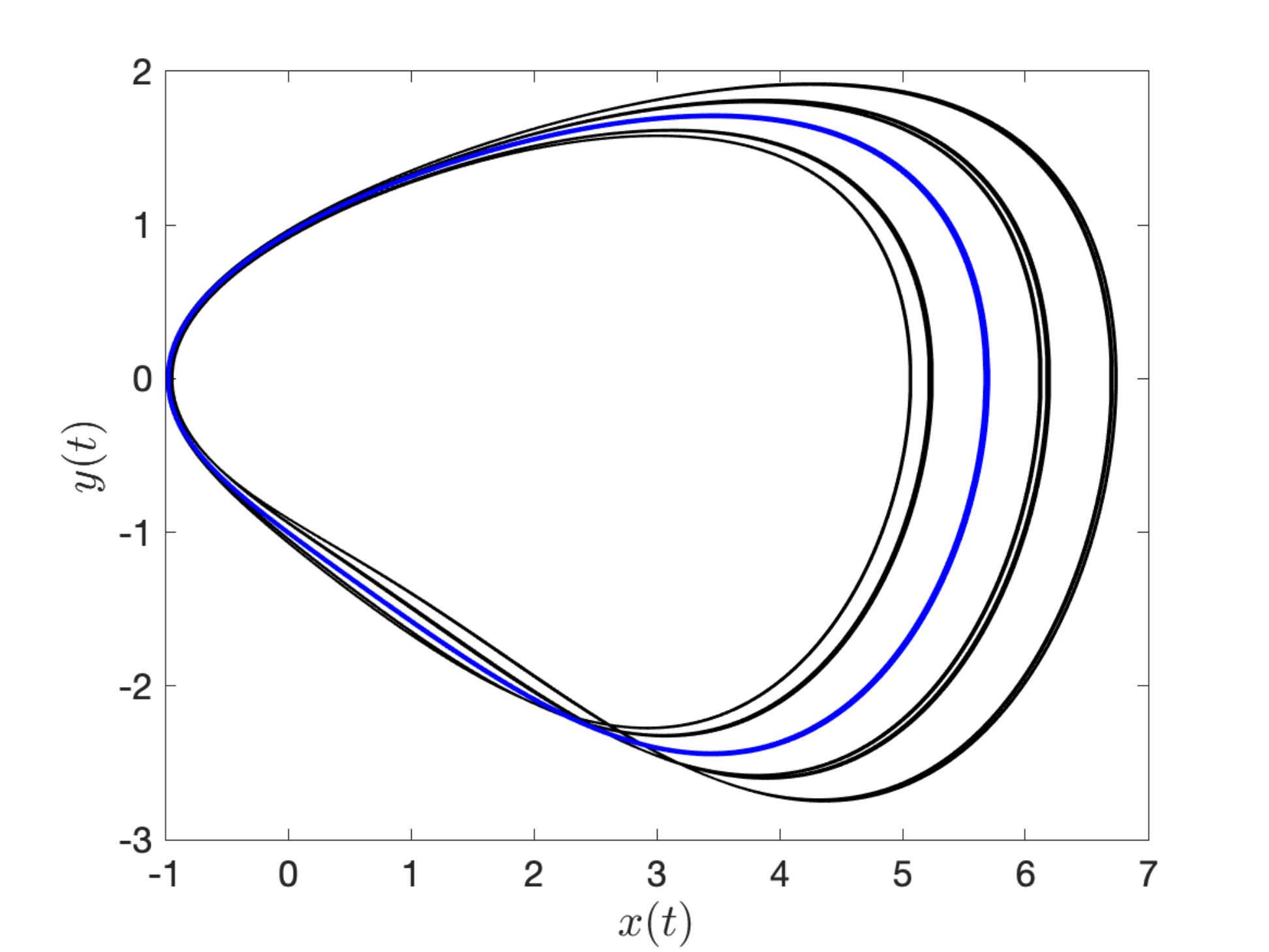} & 
\includegraphics[width = 0.35\textwidth]{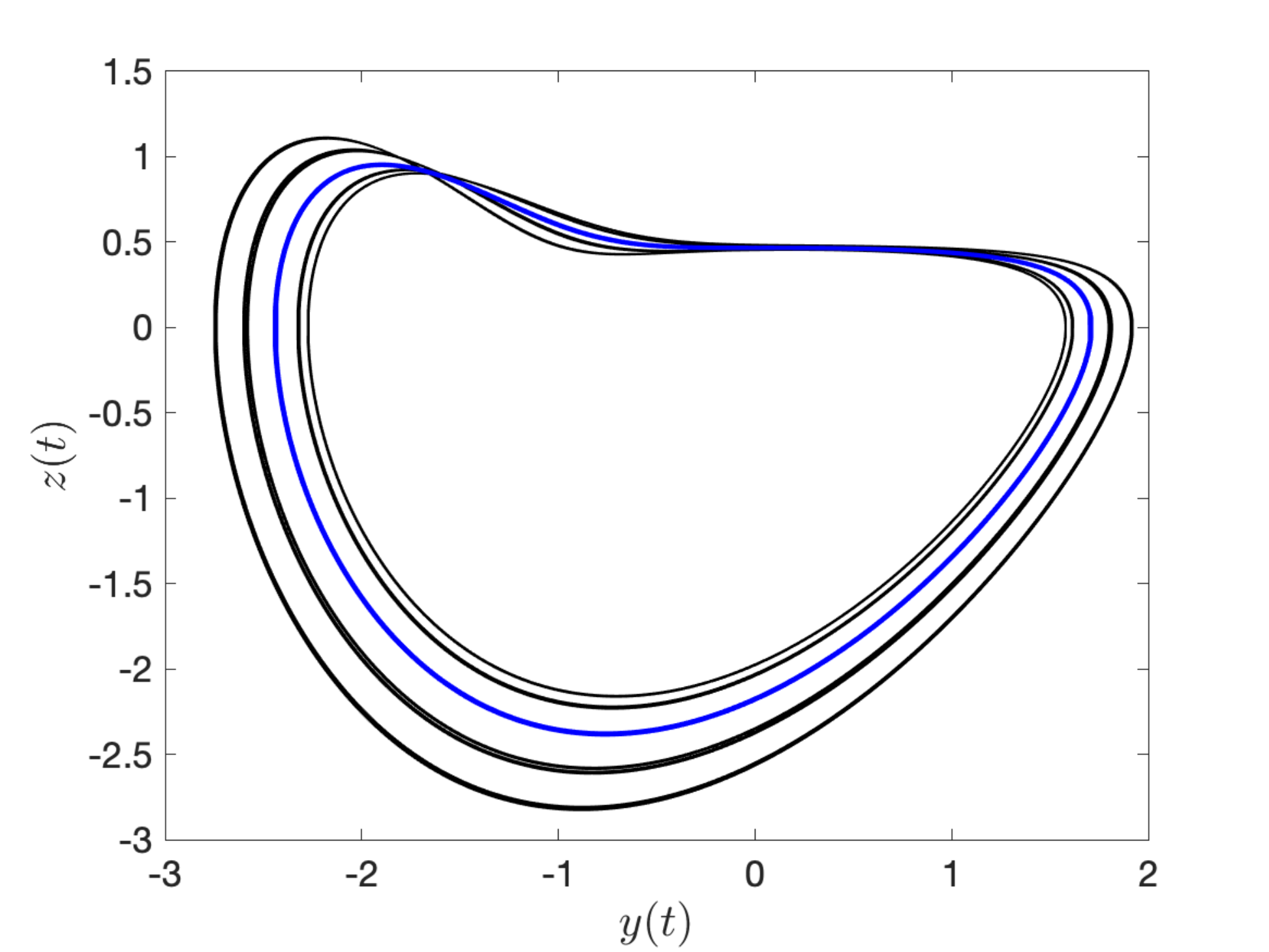} & \begin{tabular}{c} $\bar{\mu} = 2.06$ \\ $(\bar{x},\bar{z}) = (5.5228,-2.1877)$ \\ $\Kv = [0.19481\ 0.19173]$ \\ $\eta = 0.3$ \end{tabular} \\
\includegraphics[width = 0.35\textwidth]{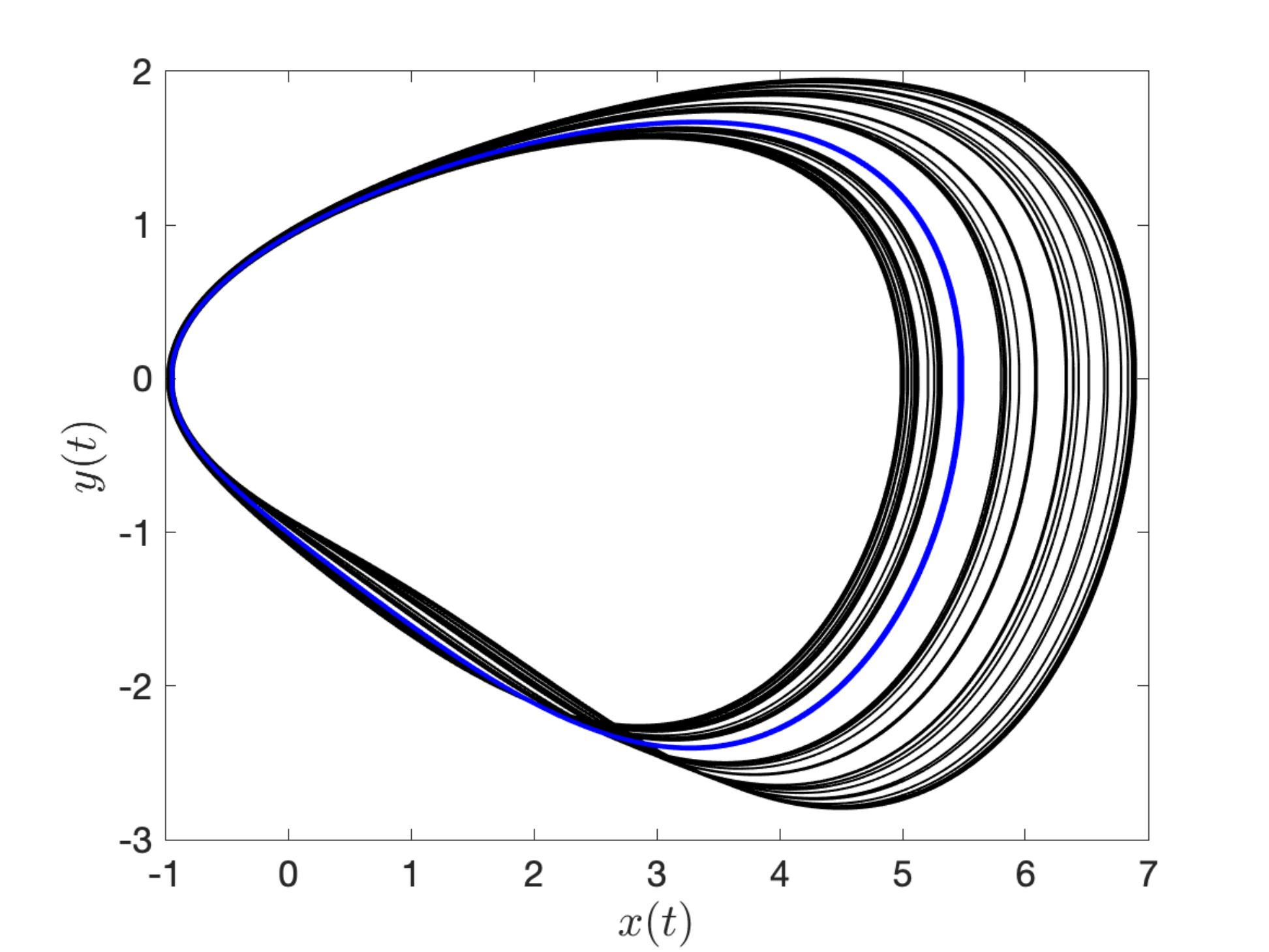} & 
\includegraphics[width = 0.35\textwidth]{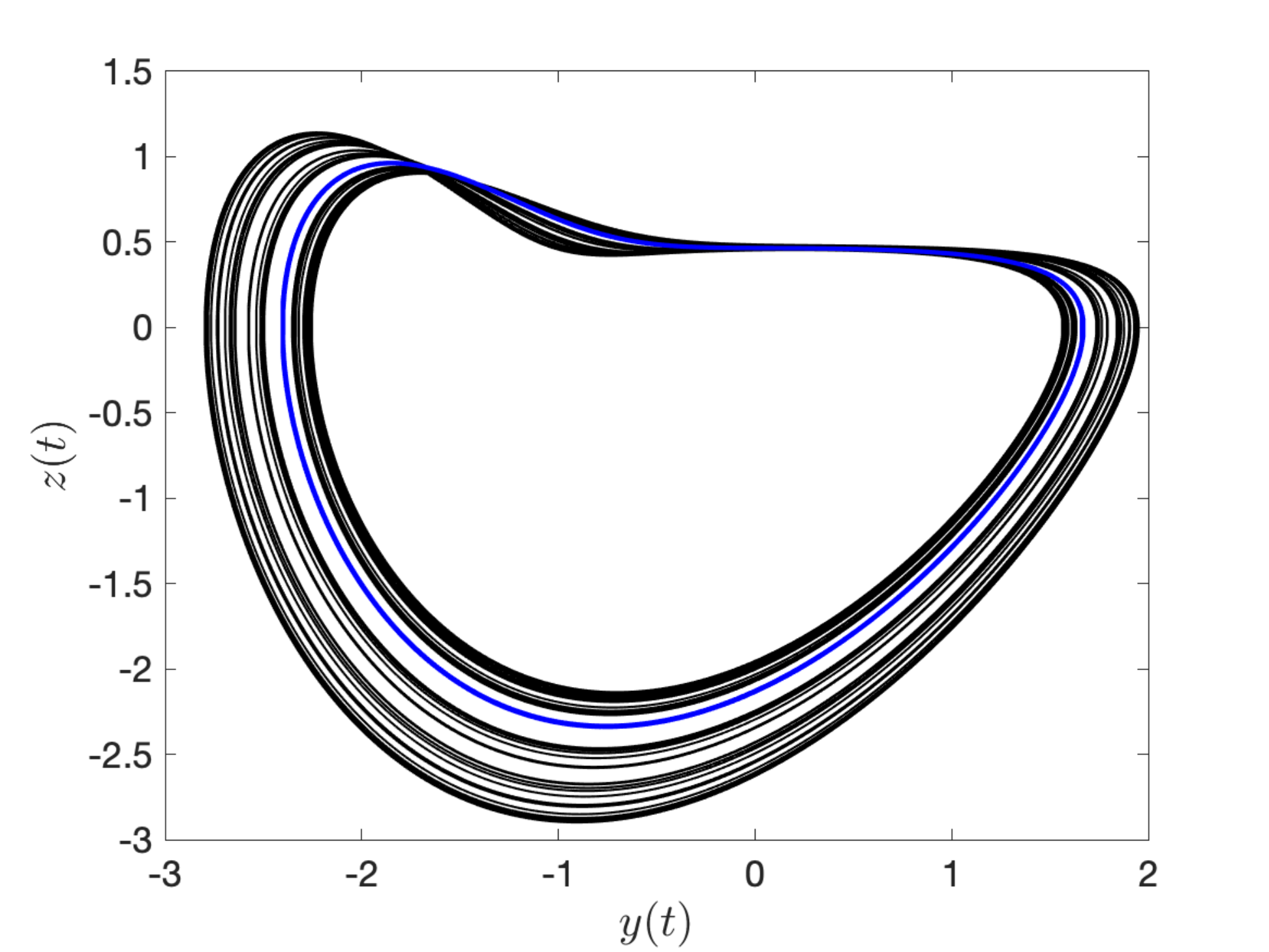} & \begin{tabular}{c} $\bar{\mu} = 2.05$ \\ $(\bar{x},\bar{z}) = (5.4803,-2.1189)$ \\ $\Kv = [0.20580\ 0.18280]$ \\ $\eta = 0.1$ \end{tabular} \\
\end{tabular}
\caption{Stabilized period 1 orbits (blue) of the Sprott system \eqref{Sprott2} along with the attractor (black) for the given parameter value. On the left we plot $y(t)$ against $x(t)$ and in the middle we plot $z(t)$ against $y(t)$. On the right we provide the focal parameter value $\bar{\mu}$, location of the fixed point in the discovered Poincar\'e map $(\bar{x},\bar{z})$ in the section $y = 0$, the control matrix $\Kv\in\R^{1\times2}$, and the threshold parameter $\eta$.}
\label{fig:Sprott}
\end{table}


\section{Closed Loop Parameter-Independent Control}\label{sec:ParamInd}

Throughout this manuscript we have focused on parameter-dependent systems and shown how data-driven discovery of Poincar\'e sections can be combined with pole-placement techniques to stabilize UPOs via parameter manipulation. Of course, many systems do not have controllable external parameters, and so here we briefly discuss the application of our stabilization algorithm using state-dependent closed loop feedback based only on the state of the system at each iteration in the Poincar\'e section. We also present an application to orbital station-keeping using energy-efficient thruster burns. 

In principle, one may consider a mapping in the Poincar\'e section of the form 
\begin{equation}\label{MapInd}
	\xv_{n+1} = \Fv(\xv_n,\uv_n),
\end{equation} 
where $\uv_n \in \R^p$ is the control variable to be applied at each map iteration. Hence, we have essentially replaced the parameter dependence of \eqref{Map} with the external control $\uv_n$. One way in which we could obtain a mapping of the form \eqref{MapInd} would be to employ the methods of Kaiser et al.~\cite{Kaiser2} which augments the system identification method of Section~\ref{sec:SINDy} to include external controls. Another way would be to have a parameter-independent Poincar\'e section mapping of the form 
\begin{equation}
	\xv_{n+1} = \Fv(\xv_n),
\end{equation}
and to apply control to vary the state of $\xv_n$ at each iteration. This amounts to iterating 
\begin{equation}\label{MapKick}
	\xv_{n+1} = \Fv(\xv_n + \Rv\uv_n),
\end{equation}
where $\Rv\in\R^{d\times p}$ is a matrix which represents potential physical restrictions on adding an impulse to the state of the system at each iteration $n$. It will be mappings of the form \eqref{MapKick} that we restrict our attention to in what follows. 

With a mapping of the form \eqref{MapKick} the process of obtaining sufficiently small kicks to the system through the choice of $\uv_n$ to stabilize a given fixed point or cyclic solution is completely analogous to the work of Section~\ref{sec:LMI}. Let us illustrate this process by considering a fixed point $\bar{\xv}$ that we wish to stabilize. Taking again a threshold parameter $\eta > 0$, we may work to find a matrix $\Kv\in\R^{p\times d}$ so that
\begin{equation}
	\uv_n = \begin{cases}
		\Kv(\xv_n - \bar{\xv}) & \|\xv_n - \bar{\xv}\|\leq \eta \\
		0 & \|\xv_n - \bar{\xv}\| > \eta
	\end{cases}	
\end{equation}   
so that the control is only applied when we are sufficiently close to the fixed point $\bar{\xv}$. Linearizing \eqref{MapKick} in a neighbourhood of $(\xv_n,\uv_n) = (\bar{\xv},0)$ gives
\begin{equation}
	\xv_{n+1} - \bar{\xv} \approx \Av(\xv_n - \bar{\xv}) + \Av\Rv\Kv(\xv_n - \bar{\xv}) = (\Av + \Av\Rv\Kv)(\xv_n - \bar{\xv}),
\end{equation} 
where again we use the notation $\Av = \Fv_x(\bar{\xv})$. From here the procedure of finding an appropriate matrix $K$ that stabilizes $\bar{\xv}$ can be obtained by implementing the LMI method of Section~\ref{subsec:FixedPt} with $\Av$ as given and $\Bv = \Av\Rv$. The extension to cyclic points follows in a nearly identical manner. 

Let us illustrate this extension of our method with a restricted 3-body problem for which two of the masses, $m_1$ and $m_2$, are much larger than the third mass, $m_3$. We consider the 2-body problem formed by $m_1$ and $M_2$ in an inertial frame of reference, restricted to their orbital plane, with $m_1$ fixed at the origin and $m_2$ fixed at $(1,0)$. Then, the governing equations of the motion of $m_3$ under the gravitational pull of the two heavy masses can be written as the four-dimensional ODE
\begin{equation}\label{3Body}
	\begin{split}
		\dot{x} &= z \\
		\dot{y} &= w \\
		\dot{z} &= 2w + x  - \frac{\mu (x - 1)}{((x-1)^2 + y^2)^{3/2}} - \frac{x}{(x^2 + y^2)^{3/2}} \\
		\dot{w} &= -2z + y  - \frac{\mu y}{((x-1)^2 + y^2)^{3/2}} - \frac{y}{(x^2 + y^2)^{3/2}}.		
	\end{split}
\end{equation}
The quantity $\mu \in (0,1)$ measures the relative masses of the smaller to large mass, i.e. $\mu = m_2/m_1$. As a specific application of our method, we will consider $m_1$ to be the Earth and $m_2$ the Moon, with $m_3$ representing a relatively massless satellite. Since the Moon is approximately $1.2\%$ of the Earth's mass, we will fix $\mu = 0.012$ in what follows.  

System~\eqref{3Body} has five equilibria, denoted $L_1,\dots,L_5 \in \mathbb{R}^4$, representing the Lagrange points of the system. These points are given to five significant digits by
\begin{equation}
	\begin{split}
		L_1 &= (0.85006,0,0,0) \\
		L_2 &= (1.1667,0,0,0) \\
		L_3 &= (-1.0010,0,0) \\
		L_4 &= (0.48814,0.86603,0,0) \\
		L_5 &= (0.48814,-0.86603,0,0).
	\end{split}
\end{equation}
The equilibria $L_1,L_2,L_3$ are saddles, and are therefore unstable, while $L_4$ and $L_5$ are stable. We will assume the satellite can be controlled externally via slight thrusts, and therefore it will be our goal to have the satellite sit at one of the unstable Lagrange points for an arbitrarily long time by applying appropriate thrusts at discrete times. This process is referred to as {\em orbital station-keeping} and has the effect that we can conserve fuel since the thrusts are only applied as slight kicks to the system which in our case will be proportional to the distance the satellite is from the Lagrange points. Under these assumptions, the desired Poincar\'e section is obtained by tracking the solutions of \eqref{3Body} at the times $t_n = nT$, where $T > 0$ represents the time between when thrusts are applied. Denoting $\xv_n = (x(t_n),y(t_n),z(t_n),w(t_n))$ to be the solution of \eqref{3Body} at time $t = t_n$, $n \geq 0$, we will use the data-driven discovery method of Section~\ref{sec:SINDy} to obtain a mapping $\Fv:\R^4\to\R^4$ that approximately gives $\xv_{n+1} = \Fv(\xv_n)$. We consider the control $\uv_n\in\R^2$ to represent thrusts in the $x$-direction with its first component and thrusts in the $y$-direction with its second. We can only add our control to the third and fourth components of $\xv$ because these velocity components are all that the impulsive thrusts can effect, so the matrix $\Rv\in\R^{4\times2}$ in \eqref{MapKick} is given by
\begin{equation}
	\Rv = \begin{bmatrix}
		 0 & 0 \\ 0 & 0 \\ 1 & 0 \\ 0 & 1
	\end{bmatrix}.
\end{equation}   
We seek an appropriate control matrix $\Kv\in\R^{2\times 4}$ that stabilizes one of the unstable Lagrange points. 

We begin by focusing on the Lagrange point $L_2$. We take $T = 0.5$, meaning that thrusts will be applied every $0.5$ time units, and generate training data from \eqref{3Body} exclusively in a neighbourhood of the desired equilibrium since we are only interested in the linearized dynamics near this point. One initial condition used for the training data is exactly at the equilibrium for the mapping to learn that this is a fixed point, while another 18 initial conditions are generated by initial conditions off this equilibrium but very close so that the saddle structure can be observed in the training data. From the discovered mapping we find that the linearized dynamics near the fixed point $L_2$ are governed by the matrix 
\begin{equation}\label{L2Matrix}
	\Av = \begin{bmatrix}
		 2.0478 & -0.093413 & 0.53058 & 0.21987 \\
  		-0.34021 & 0.76162 & -0.22646 & 0.39881 \\ 
   		4.4265 & -0.55106 & 1.3841 & 0.98674 \\ 
  		-1.9588 & -0.81595 & -0.99632 & 0.30130
	\end{bmatrix}.
\end{equation}
The matrix \eqref{L2Matrix} has eigenvalues $0.33348, 2.9612, 0.60008 \pm 0.80248{\rm i}$, so that one lies outside the unit circle of the complex plane, one inside, and two approximately on it. Thus, $L_2$ is unstable in the obtained coarse-grained mapping and just like in \eqref{3Body} it (approximately) has one unstable direction, one stable direction, and a two-dimensional centre manifold. Using the LMI procedure of Section~\ref{sec:LMI} we obtain the control matrix $\Kv \in \R^{2\times 4}$ given by
\begin{equation}\label{3BodyK}
	\Kv = \begin{bmatrix}
		-2.6242 & 0.24424 & -0.95380 & -0.051993 \\
		-1.0448 & 0.061818 & -0.032256 & -0.75239
	\end{bmatrix}
\end{equation}
which gives that $\Av + \Av\Rv\Kv$ now has eigenvalues $0.048357,0.47768, 0.22383\pm 0.21157{\rm i}$. Hence, $\Av+\Av\Rv\Kv$ is stable, as desired.

\begin{figure} 
\center
\includegraphics[width = 0.49\textwidth]{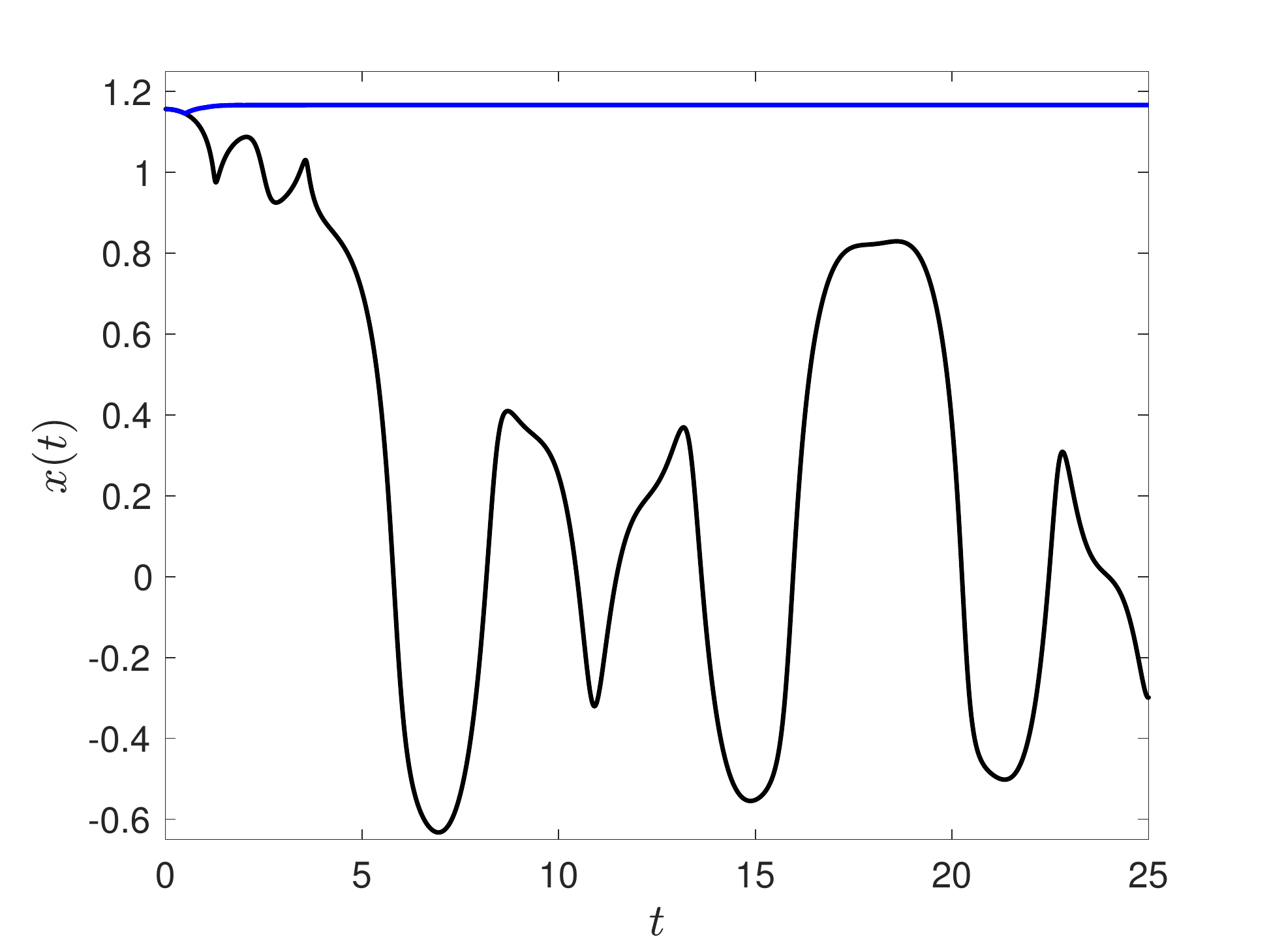} \ 
\includegraphics[width = 0.49\textwidth]{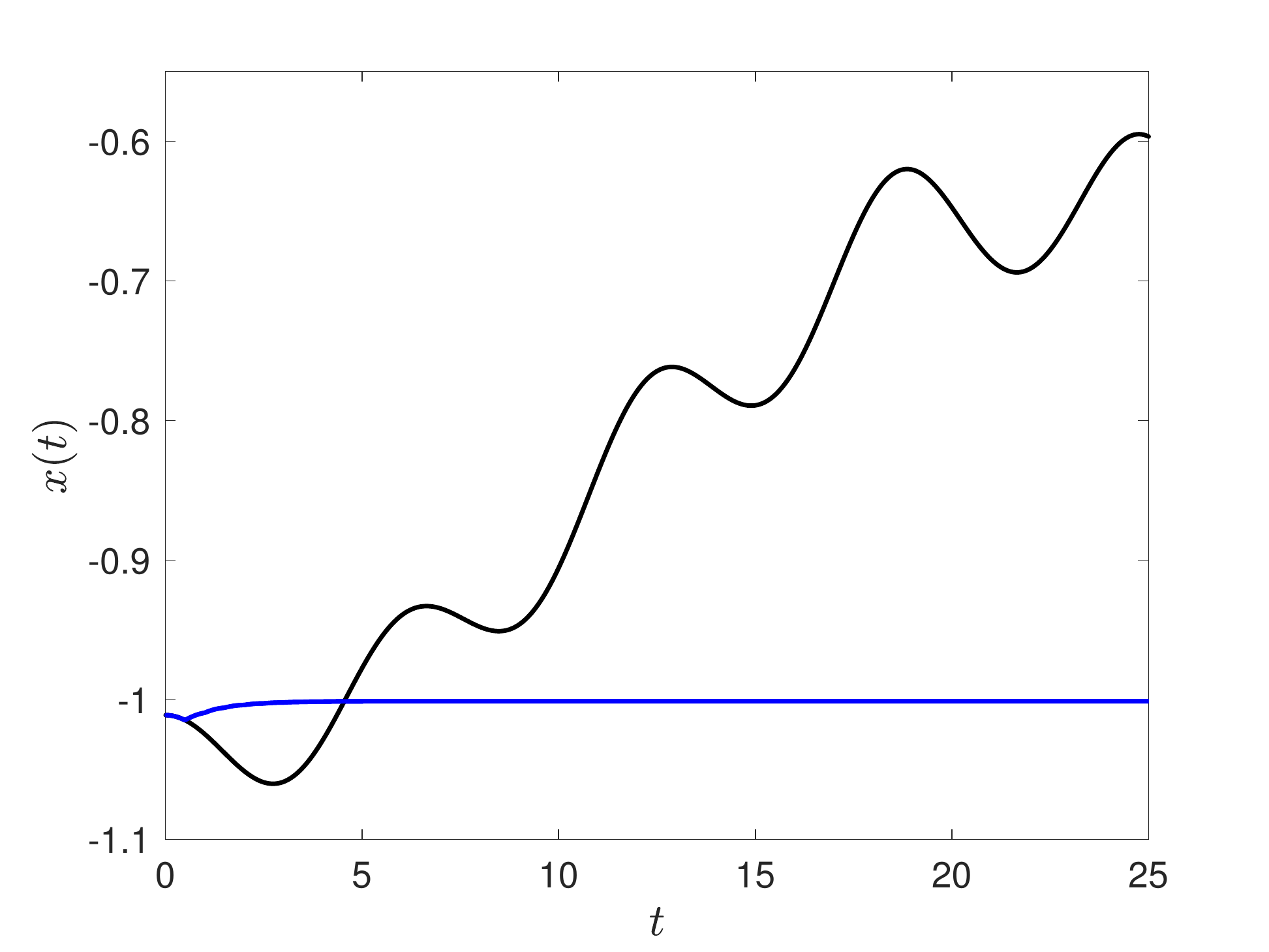}
\caption{Left: The controlled trajectory (blue) remains near the Lagrange point $L_2$ of \eqref{3Body} while the uncontrolled trajectory with the same initial conditions (black) almost immediately leaves a neighbourhood of this saddle point. Right: The same as the left image except for the Lagrange point $L_3$.}
\label{fig:3Body}
\end{figure}

This application of our method bears a slight resemblance to the stabilization of the isolated periodic orbit in \S~\ref{subsec:Hopf} in that trajectories exponentially separate from the element which one wishes to stabilize, potentially not to return. Therefore, to remain close to the Lagrange point $L_2$ of \eqref{3Body} we are required to start very close to it. In Figure~\ref{fig:3Body} we present a controlled trajectory using the matrix \eqref{3BodyK} and initial condition $(x(0),y(0),z(0),w(0)) = L_2 - (0.01,0,0,0)$ compared with an uncontrolled trajectory with the same initial condition. The threshold parameter is set to $\eta = 0.1$ since we can see that for small $t\geq 0$ the $x(t)$ trajectory quickly diverges from the Lagrange point, thus necessitating a sufficiently large control threshold to sufficiently thrust the system back towards the equilibrium. Finally, we comment that we attempted to stabilize the trajectory with a larger time between thrusts, i.e. $T > 0.5$, but with $T = 0.75$ and $T = 1$ one is required to start extremely close to the equilibrium in order to control the trajectory. The reason for this is trajectories diverge quickly from the saddle point $L_2$ and therefore the longer one waits to apply the first thrust, the farther from this equilibrium one can travel. 

The process of orbital station-keeping near the Lagrange point $L_3$ is nearly identical to that of $L_2$ and so we only plot the result of a controlled versus uncontrolled trajectory in Figure~\ref{fig:3Body} with initial condition $(x(0),y(0),z(0),w(0)) = L_3 - (0.01,0,0,0)$. We have found that applying the procedure to keep an orbit close to $L_1$ will fail if $T = 0.5$, thus requiring thrusts that are spaced closer together in time. The reason for this is that the rate of separation from the Lagrange point $L_1$ is significantly faster than that of the Lagrange points $L_2$ and $L_3$. We find success when taking $T = 0.25$, thus doubling the frequency of thrusts from the stabilization procedure for $L_2$ and $L_3$. With these thrusts at every $0.25$ time units we are able to stabilize the satellites orbit near $L_1$ for arbitrarily long times, in the same way as $L_2$ and $L_3$. The reader is directed to the accompanying code at {\bf GitHub/jbramburger/Stabilizing\_UPOs} to observe the success of our stabilization procedure applied to $L_1$ since a figure is not provided here in an effort to avoid redundancy.


\section{Discussion}\label{sec:Discussion}

In this work we have presented a method for stabilizing UPOs of ordinary differential equations. Our approach uses the recently developed SINDy method~\cite{SINDy,Bramburger} to discover a parsimonious representation of a Poincar\'e mapping which can be used to both find and analyze the stability of periodic orbits. Once an unstable fixed point or cyclic orbit has been identified in the Poincar\'e mapping, we may use the pole-placement method of Romeiras et al.~\cite{Romeiras} to apply slight parameter adjustments each time a trajectory intersects the Poincar\'e section to stabilize a UPO. Furthermore, we have demonstrated how to automate the process of obtaining the appropriate controls to stabilize a fixed point using Parrilo's LMI framework~\cite{Parrilo}, while also extending these ideas to obtain controls for cyclic orbits. 

We applied our method to a number of systems, each presenting a slightly different perspective on the application of the control algorithm. Our work began with the chaotic H\'enon map where we focussed exclusively on the LMI procedure for obtaining appropriate controls to stabilize a fixed point, 2-cycle, and 4-cycle. We then moved to a simple planar ODE that had an isolated periodic orbit which we were able to stabilize by first discovering a sparse representation of a Poincar\'e mapping and then determining the appropriate parameter manipulations to stabilize the orbit. In \S~\ref{subsec:Rossler} we stabilized a number of periodic orbits in the R\"ossler system which has the advantage of having a one-dimensional Poincar\'e section. We demonstrated how following period-doubling bifurcations with the training data can be used to help the discovered SINDy mapping infer the existence of UPOs, while also showing that the discovered mapping can do a good job of finding these UPOs even if the training data does not contain snapshots of its destabilizing bifurcation. Both the R\"ossler system and Sprott's chaotic jerk system provided examples of the performance of the method in fully chaotic regimes, where we were able to stabilize both period 1 and 2 orbits. Finally, in Section~\ref{sec:ParamInd} we demonstrated the extension of this method to systems that do not have explicit parameter dependence but can be controlled externally based on instantaneous system measurements. Our methods were applied to a classic problem in control theory: having a satellite rest at a Lagrange point in a restricted 3-body problem via impulsive thrusts applied at evenly spaced times.         

There are a number of ways to improve this method moving forward. As we saw in many of our examples, knowing the location of the UPO that one wishes to stabilize is advantageous since this information can be added to the training data and therefore reflected in the discovered mapping. Hence, upon gathering section data from a chaotic trajectory of an ODE, we may apply the work of So et al.~\cite{So} to extract UPOs from this section data. This would be especially useful when the section data is gathered from a real-world system where regeneration of the data with multiple initial conditions is not necessarily practical. Beyond this, new advances using dynamic mode decomposition~\cite{Page} and sum-of-squares relaxations of differential inequalities~\cite{Lakshmi} have yielded methods of obtaining UPOs of dynamical systems, thus providing another method by which we could build information about UPOs into the training data. Sum-of-squares methods may also be useful for determining the size of the basin of attraction for a stabilized orbit of the Poincar\'e mapping~\cite{Sidorov}, thus providing upper bounds on the size of the threshold parameter $\eta > 0$. 

Moving forward it is desirable to extend these methods to infinite-dimensional differential equations, particularly spatially extended systems. One method is to follow in a similar manner to the pioneering work of Lorenz~\cite{Lorenz} by projecting the dynamics of a partial differential equation onto a finite collection of elements of a Fourier or Galerkin basis~\cite{Noack2003jfm}. This would result in an ordinary differential equation for the coefficients of the basis elements that can be controlled with our methods. The drawback to this method is that unless the dynamics of the full system in question are genuinely finite-dimensional, it is likely that only controlling a finite number of basis elements will not give way to control of the infinite-dimensional dynamical system. Another limitation here would be that as the number of degrees of freedom in a system increases, the number of variables for discovering polynomial mappings increases exponentially. Hence, it appears that controlling UPOs in spatio-temporal systems whose dynamics are not completely described by finite-dimensional dynamics require nontrivial extensions to this method that would potentially utilize a dimensionality reduction component. We hope to report on this in a follow-up investigation.  
Finally, there are a number of compelling applications that may benefit from this control approach, including space mission design~\cite{koon2000heteroclinic,gomez2004connecting,Dellnitz:2005} and fluid flow control~\cite{Budanur,Cvitanovic2,Fazendeiro,Franceschini,Lucas,petrov1996nonlinear,eckhardt2007turbulence,Yalniz,suri2020capturing,graham2020exact}.  These topics are the subject of ongoing work, although SINDy has recently been applied to several fluid flows with promising results~\cite{Loiseau2017jfm,Loiseau2018jfm,loiseau2020data}.   

\section*{Acknowledgements}
JNK acknowledge support from the Air Force Office of Scientific Research (AFOSR) (FA9550-17-1-0200).
SLB acknowledges funding from the Army Research Office (ARO W911NF-19-1-0045).  

\setlength{\bibsep}{2.6pt plus 1ex}
\begin{spacing}{.01}
	\small
\bibliographystyle{unsrt}
\bibliography{Stabilizing_UPOs}
\end{spacing}
\end{document}